\input amsppt.sty
\def\RefKeys{ B BMM BS Ca Fi Fu G1 G2 G3 G4 Gas GK H HL HMS KT1 KT2 L Li Ma
MO LeT
LT1 LT2 Pa R S T1 T2 T3 T4 V }
        \widestnumber\key{BMM}

\define \OO{\Cal O}
\define \II{\Cal I}

\define \C{\Bbb C}

\define \PP{\Bbb P}
\define \V{{\roman V}}

\define \Bl{\roman{ Bl}}

\define \Cyc {\roman{ Cyc}}
\define \mult {\roman{mult}}
\define \dist{\roman{dist}}

\define \sub{\subseteq}

\define \inv{^{-1}}
\define \pd{\partial}
\define \pdz#1{\pd_{z_#1}}

\define \pda#1{\pd_{a_{#1}}}




\def\TheMagstep{\magstep1}      
 \def\PaperSize{letter}         

\def\TRUE{TRUE}
\ifx\DoublepageOutput\TRUE \def\TheMagstep{\magstep0} \fi
\mag=\TheMagstep
\TagsOnRight
\vsize48pc
\abovedisplayskip 6pt plus6pt minus0.25pt
\belowdisplayskip=\abovedisplayskip
\abovedisplayshortskip=0mm
\belowdisplayshortskip=2mm
\def\centertext
 {\hoffset=\pgwidth \advance\hoffset-\hsize
  \advance\hoffset-2truein \divide\hoffset by 2
  \voffset=\pgheight \advance\voffset-\vsize
  \advance\voffset-2truein \divide\voffset by 2
 }
\newdimen\pgwidth	\newdimen\pgheight
\def\letter{letter}	
\ifx\PaperSize\letter
	\pgwidth=8.5truein \pgheight=11truein \centertext \fi
\ifx\PaperSize\AFour
	\pgwidth=210truemm \pgheight=297truemm \centertext \fi

\def\atref#1{\ref\key\csname#1\endcsname}
 \newcount\refno \refno=0
 \def\MakeKey{\advance\refno by 1 \expandafter\xdef
 	\csname\TheKey\endcsname{{%
	\ifx\UseNumericalRefKeys\TRUE
		\number\refno
	\else\TheKey \fi}}\NextKey}
 \def\NextKey#1 {\def\TheKey{#1}\ifx\TheKey\NoKey\let\next\relax
  \else\let\next\MakeKey \fi \next}
 \def\NoKey{*!*}
 \expandafter\NextKey \RefKeys *!*
 \ifx\UseNumericalRefKeys\TRUE \widestnumber\no{\number\refno}\fi

 \newskip\sectskipamount \sectskipamount=0pt plus30pt
 \def\sectionhead#1 #2\par{\def\sectno{#1}\def\sectname{#2}%
   \vskip\sectskipamount\penalty-250\vskip-\sectskipamount
   \bigskip
   \centerline{\smc \number\sectno.\enspace\sectname}\nobreak
   \medskip
   \message{\number\sectno. \sectname }%
}

 \def\proclaimm#1#2 {\medbreak\noindent {\bf(\sectno.#2) #1.}\quad
	\begingroup\it}
 \def\subh#1#2{\medbreak\noindent {\bf(\sectno.#2) #1.}\quad}
\def\art#1 #2\par{\subh{\rm({\it#2\unskip\/})}{#1}}
\def\rmk{\subh{Remark}}

\def\stp{\subh{Setup}}
\def\expl{\subh{Example}}
\def\prop{\proclaimm{Proposition}}
	
\def\thm{\proclaimm{Theorem}}
	\def\thmx #1{\proclaimm{Theorem {\rm (#1)}}}
\def\cor{\proclaimm{Corollary}}
\def\lem{\proclaimm{Lemma}}
	
\def\pf{\endgroup\medbreak\noindent {\bf Proof.}\quad}
\def\demobox{\vbox{\hrule\hbox{\vrule\kern.5ex
	\vbox{\kern1.2ex}\vrule}\hrule}}
\def\enddemo{{\unskip\nobreak\hfil\penalty50
  \hskip1em\hbox{}\nobreak\hfil\demobox
  \parfillskip=0pt \finalhyphendemerits=0 \par}}

\def\tgs{\tag\sectno.}

\centertext
\topmatter
\title Segre Numbers and Hypersurface Singularities
\endtitle

\author Terence Gaffney and Robert Gassler  \endauthor


\thanks{The first author was partially supported by N.S.F.~grant
9403708-DMS.}
\endthanks

\endsubjclass

\abstract
We define the Segre numbers of an ideal as a generalization of the  multiplicity
of an ideal of finite colength. We prove generalizations of various theorems
involving the multiplicity of an ideal such as a  principle of specialization of
integral dependence, the Rees-B\"oger theorem, and the  formula for the
multiplicity of the product of two ideals. These results are applied to the
study
of various equisingularity conditions, such as Verdier's condition W, and
conditions $A_f$ and $W_f$.
\endabstract
\endtopmatter

\define \intro{0}
\def \class{1}
\define \segre{2}
\define\form{3}
\define \sid{4}
\define\af{5}
\define\wf{6}

\define\th#1{\tilde h_{#1}}
\define\tD#1{\tilde D_{#1}}
\define\tb#1{\tilde b_{#1}}

\centerline{alg-geom/9611002}\medskip


\sectionhead {\intro} Introduction


If an ideal $I$ in a Noetherian local ring $A$ has finite colength,
 then the multiplicity of the ideal is  a fundamental invariant of the ideal
with many applications in geometry and algebra.  Pierre Samuel \cite{\S}
 used it to define the intersection multiplicity of two algebraic sets.
David Rees
\cite{\R} linked the multiplicity of $I$ to its integral closure $\bar I$.
Bernard Teissier \cite{T1} used it to study the equisingularity of families of
hypersurfaces with isolated singularities.

If an ideal does not have finite colength then the multiplicity is not
defined.   In this paper we study ideals of non-finite colength, and we define
a set of invariants associated to $I$,
 the Segre numbers of $I$. This set of invariants then has similar properties
to the multiplicity. If the multiplicity of $I$ is defined, it can be realized
as the degree of the exceptional
 divisor in the blowup by $I$. The Segre numbers are constructed by first
forming cycles by intersecting the exceptional  divisor in the blowup by $I$
with generic hyperplanes, then  pushing the cycles down to the base,  and
finally taking multiplicities of these cycles.

The Segre numbers provide a link between
$I$ and $\bar I$ just as in the case of ideals of finite colength;
 see Corollary (\sid.9). They are
lexicographically upper semi--continuous, as observed in (\sid.5).
They are also useful in studying the equisingularity of families of
hypersurfaces with non-isolated singularities.

Suppose we are given a family of hypersurfaces $X=\{X(t)\}$ parameterized by
$(\C^p,0)$ and
 defined by  $f:(\C^{n+1+p},0)\to(\C,0)$ a function, which vanishes on
$0\times(\C^p,0).$ The goal is to find invariants which depend only on the
members of the family $X(t)=V(f_t)$ whose constancy ensures that the given
equisingularity condition being studied holds for the family.  The blowup of
$(\C^{n+1+p},0)$ along the Jacobian ideal $J(f)$ and the blowup of
$X$ along the ideal $J(f)'$ induced by the Jacobian ideal are useful
starting points for many equisingularity conditions.  These spaces are the
relative conormal of $f$ and the conormal of
$X.$ The latter consist of  the closure of
the pairs $(x,H)$ where $x$ is a smooth point of $X$ and $H$ is the tangent
hyperplane to $X$  at such a point, while the relative conormal is the closure
of the pairs $(x,H)$, where $x$ is  a point $f$ is a submersion, and $H$
is the tangent hyperplane to the fiber of $f$ through $x$.  The exceptional
divisors of these blowups record the behavior of the limiting tangent
hyperplanes to $X$, resp.  the smooth fibers of $f$. The equisingularity
conditions that we have in mind, Whitney conditions and the
$W_f$ condition, are defined in terms of the behavior of these limiting tangent
planes.

The Segre numbers of $J(f_t)$ and the ideal $J(f_t)'$ induced by it in
$\OO_{X(t),0}$ give us invariants which depend only on $X(t)$ and which
describe the exceptional divisors of $\Bl_{J(f_t)}(\C^{n+1},0)$ and
$\Bl_{J(f_t)'}X.$  To show that an equisingularity condition holds, it is
necessary to pass to controlling the exceptional divisors of
$\Bl_{J(f)}(\C^{n+1+p},0)$ and
$\Bl_{J(f)'}X$.  More generally, for an an arbitrary family $X$ and an ideal
$I$ in $\OO_{X,0},$ we want to control the exceptional divisor of
$\Bl_{I}X$ by conditions imposed on the exceptional divisors of
$\Bl_{I(t)}X(t).$  This is done in Theorem (\sid.6) which is a
generalization of the principle of
specialization of integral dependence due to Teissier. This theorem has as a
corollary (\sid.9) a generalization of the theorem of Rees mentioned
earlier. B\"oger's theorem
 then follows as an easy application (\sid.10).

It turns out that the Whitney conditions are controlled by the Segre numbers
of
$mJ(f_t)$, where $m$ is the maximal ideal in $\OO_{\C^{n+1},0}.$  It is useful
to relate these Segre numbers to the Segre numbers of $J(f_t)$ since these
are more closely linked to the  geometry of $X(t)$. In fact, the Segre numbers
of  $J(f_t)$ are just the L\^e numbers of David Massey \cite{\Ma},
 and their alternating sum is the Euler characteristic of the Milnor fiber of
$f$. This is done in greater generality in Theorem (\form.5)   which
relates the Segre
numbers of the product of an ideal $I$ with an $m$--primary ideal $J$ to the
Segre numbers of $I$ and other invariants of the pair
$I$ and $J$.  This extends a result of Teissier \cite{T1}. The theorem of
Teissier relates the multiplicity of the product of two ideals of finite
colength with their mixed multiplicies. In our theorem the mixed
multiplicities are seen to be the polar multiplicities of the ideal $I$
computed with respect to the ideal
$J$.  In the case where $I=J(f_t)$  and $J=m$ these polar
multiplicities are just the relative polar multiplicities of $f_t$.

All of these results come together in the main results of section \wf. There we
show that the  constancy of the relative polar multiplicities of $f_t$ and of
the L\^e numbers of $f_t$ implies that
$f$ satisfies the $W_f$ condition and that the smooth strata of $X$, its
singular set $\Sigma(X)$ and the components of the singular set $\Sigma
(\Sigma(X))$ of the singular set which are of codimension 1 in  $\Sigma(X)$ are
all Whitney regular over the parameter stratum. We also show that if $X$ has a
Whitney stratification which includes
$0\times(\C^p,0)$ as a stratum,
 then the above numbers are constant. This gives as a corollary a theorem of
Adam Parusinski \cite{\Pa} which says that the existence of such a Whitney
stratification for $X$ implies the $W_f$ condition  for $f$. This proof
fulfills  Parusinski's prediction that his result could be proved using  the
L\^e numbers.

These results are possible because of a surprising phenomenen. It turns out
that the constancy of
the relative polar multiplicities of $f_t$ and of
the L\^e numbers of $f_t$ imply that the Whitney conditions hold for any
stratum $W$ in the
critical set of $f$  over the parameter stratum whose closure is the image
of a component of the
exceptional divisor of the blowup of the ambient space by the Jacobian
ideal of $f$. This means
that it is possible to show that many strata satisfy the Whitney conditions
using only these two
sets of invariants. It remains to find a good criteria for when strata
correspond to components of
the exceptional divisor.

As a further application of our results, we show that the constancy of the
relative polar multiplicities
and of the L\^e numbers in a family of 2-dimensional hypersurfaces
implies that the family is Whitney equisingular (Corollary (\wf.6)).  This
result is possible, because the failure of
smooth points of the total space to be Whitney over the smooth points of
the singular set implies the existence
of a component of the exceptional divisor of the  relative Nash blowup of
$\C^3\times\C$ of the type controlled by our numerical
conditions. This example indicates the importance of relating the
components of the exceptional divisor of this blowup to the
geometry of $X$.  We then apply
our results to the to the study of the  family of hypersurfaces given by
hyperplane slices of a three dimensional hypersurface
$X$ (Example (\wf.7)).

Finally, as an indication of the power of our approach, in (Proposition
(\wf.9)) we show that if a hypersurface satisfies the same
relations among its  L\^e numbers as a hypersurface defined by a
homogeneous polynomial, then the  smooth points of the
deformation of the hypersurface
to its tangent cone satisfy the Whitney conditions over the
parameter axis, provided the tangent cone is reduced.  If the original
space is 2-dimensional, then  by (\wf.6),  the deformation is Whitney
equisingular.

\medskip
It is a pleasure to acknowledge our debt to Bernard Teissier. His Cargese
paper \cite{T1} convinced us  to try to control equisingularity conditions by
invariants depending only on the members of the family under consideration,
and showed us the elements from the isolated case which we needed to
generalize to handle non-isolated singularities. David Massey's  work on the
L\^e numbers was a constant guide to us as we developed the properties of the
Segre numbers.  Steven Kleiman and David  Massey also supported us through
many conversations over the years that this paper developed. We also thank
Steven Kleiman for many helpful comments on earlier versions of this paper.

\medskip
The paper is organized as follows. In section \class, we recall the
basic results of the theory of  integral closure of ideals and their
application to the equisingularity theory of families of  hypersurfaces with
isolated singularities.  In section \segre, we define the Segre cycles and
polar varieties  of an ideal, and prove some of their basic properties. In
section
\form, we give a formulation of these notions using intersection theory
(\form.1), and use this to describe how the  Segre numbers change under
hyperplane section in (\form.4).  Theorem (\form.5) is the expansion formula.
In section {\sid}, we discuss the upper semi--continuity properties of the
Segre numbers, and the behavior of the polar varieties of $I(t)$ in a family.
In (\sid.6) we give a necessary and sufficient criterion for the Segre numbers
of
$I(t)$ to be constant along the parameter stratum in terms  of the behavior of
the components of the exceptional divisor of $\Bl_{I}X$ which project to the
parameter stratum, and the limiting secant behavior of the components  which
do not project to the parameter stratum. Our generalization of the principle
of  specialization of integral dependence (\sid.7) follows easily from this
result.  We then discuss the theorems of Rees and B\"oger. In section {\af}
we use our machinery to study  Thom's $a_f$ condition,
recovering a result of  Massey, while in section {\wf} we apply our results
to the Whitney conditions and the
$W_f$ condition.


\sectionhead {\class} The classical theory


\art1 Integral dependence

Integral dependence is used in local analytic
geometry to translate inequalities between analytic functions into algebra. The
basic source for this is the work \cite{LeT} of Monique Lejeune--Jalabert and
 Teissier. Another useful reference is Joseph Lipman's more algebraic survey
\cite{\Li}.

 Let $(X,0)\sub (\C^N,0)$ be a reduced analytic space germ. Let $I$ be an ideal
in the local ring $\OO_{X,0}$ of
$X$ at 0, and $f$ an element in this ring. Then $f$ is integrally dependent
on $I$ if one of the following equivalent conditions obtain:
\itemitem{(i)} There exists a positive integer $k$ and elements $a_j$ in
$I^j,$ so that $f$ satisfies the relation
$f^k+a_1f^{k-1}+\dots+a_{k-1}f+a_k=0$ in
$\OO_{X,0}.$
\itemitem{(ii)} There exists a neighborhood $U$ of 0 in $\C^N,$ a positive
real number $C,$ representatives of the space germ $X,$ the function
germ $f,$ and generators
$g_1,\dots,g_m$ of
$I$ on
$U,$ which we identify with the corresponding germs, so that for all $x$ in
$X$ the following equality obtains:
$|f(x)|\leq C \max\{|g_1(x)|,\dots,|g_m(x)|\}.$
\itemitem{(iii)} For all analytic path germs $\phi:(\C,0)\to(X,0)$ the
pull--back
$\phi^*f$ is contained in the ideal generated by $\phi^*(I)$ in the local ring
of $\C$ at 0.

As an example take $I$ to be the ideal generated by the images of
the coordinates
$z_1,\dots,z_{N-1}$ in $\OO_{X,0}$ and $f$ the image of $z_N$
in the same local ring. Then, $f$ is integrally dependent on $I$ if, and only
if, the $z_N$--axis is not contained in the tangent cone of $X$ at 0.  For if
$\phi(t):(\C,0)\to(X,0)$ is an analytic path then for suitable $k$ the limit
$\lim_{t\rightarrow 0} 1/t^{k} (z_1,\dots,z_{N})$ gives the direction of the
limiting secant line to $X$ determined by $\phi$. So none of these limits
are tangent to the
$z_N$--axis if, and only if,  the order in t of $\phi^*z_N$ is greater than
or equal to the
smallest of the orders of
$\phi^*z_i.$ This is equivalent to
$\phi^*z_N$  being contained in the ideal generated by
$\phi^*z_1,\dots,\phi^*z_{N-1}$ in the local ring of $(\C,0).$  Thus, condition
(iii) shows the claim.

\smallskip
If we consider the normalization $\bar B$ of the blowup $B$ of $X$ along the
ideal
$I$ we get another equivalent condition for integral dependence. Denote the
pull--back of the exceptional divisor $D$ of $B$ to $\bar B$ by $\bar D.$
\itemitem{(iv)} For any component $C$ of the underlying set of $\bar D,$ the
order of vanishing of the pullback of $f$ to $\bar B$ along
$C$ is no smaller than the order of the divisor $\bar D$ along $C.$

The elements $f$ in $\OO_{X,0}$ that are integrally dependent on $I$ form the
ideal $\bar I,$ the integral closure of $I.$ Often we are only interested in
the properties of the integral closure of an ideal $I;$ so we may replace $I$
by an ideal $J$ contained in $I$ with the same integral closure as $I.$ Such
an ideal $J$ is called a reduction of
$I.$ In the above example, the ideal generated by the
first $N-1$ coordinate functions in the local ring of $X$ at 0 is a reduction
of its maximal ideal if, and only if, the $z_N$--axis is not contained in the
tangent cone of $X$ at 0.

It is easy to see that  $J$ is a reduction of $I$ iff there exists a finite map
$\Bl_IX\to\Bl_JX.$ In particular, the fibres of the two blowups over 0 have the
same dimension. Samuel proved that any ideal $I$ has a reduction generated by
at most $n$ elements where $n$ is the dimension of $X$ at 0. In fact, $n$
generic linear combinations of given generators of $I$ generate a reduction.

In general, an ideal $I$ has a reduction generated by $m+1$ elements where $m$
is the dimension of the fibre of $\Bl_IX$ over 0.

\medskip
If $X$ is equidimensional, and $I$ is primary to the maximal ideal, its
integral closure is completely determined by its multiplicity $e(I).$ This
follows from the following theorem of Rees'.

\thm2 An ideal $J$ contained in $I$ is a reduction of $I$ if, and only if,
the multiplicity $e(I)$ equals the multiplicity $e(J).$\endgroup

\smallskip
Note that for any ideal $J\sub I$ the multiplicity of $J$ is not smaller than
$e(I).$

\art3 Specialization of integral dependence

In his paper
\cite{T1}  Teissier established the {\it Principle of Specialization
of Integral Dependence of Ideals:}  Consider a reduced equidimensional
family $X\to Y$ of  analytic spaces, and an ideal sheaf $\II$ on $X$
with finite co--support over $Y.$ Suppose $h$ is a section of $\OO_X$ so that
for all $t$ in a Zariski--open dense subset of $Y$ the induced section of
$\OO_{X(t)}$ on the fibre over $t$  is integrally dependent on the induced ideal
sheaf $\II\OO_{X(t)}.$ If the multiplicity $e(\II\OO_{X(t)})$ is independent
of $t$ in $Y$, then $h$ is integrally dependent on $\II.$ (The multiplicity
$e(\II\OO_{X(t)})$ is the sum of the multiplicities of the ideals induced by
$\II$ in the local rings $\OO_{X(t),x_t^i},$ where the points
$x_t^1,\dots,x_t^{k(t)}$ form the fibre of the support of $\OO_X/\II$
over $t.$)

\smallskip
We give a brief sketch of the proof the principle as contained in
\cite{T3, Appendix 1}. For simplicity assume that
$\OO_X/\II$ is supported on the image of a section $t\mapsto x_t$ of $X$
over
$Y.$ Denote the fibre dimension of $X\to Y$ by $n.$ By Samuel's theorem, we can
find $n$ elements $g_1,\dots,g_n$ in the stalk of $\II$ at $x_ 0$ the image of
which in
$\OO_{X(0),x_0}$ generate a reduction of the ideal induced by $\II.$ We
identify all germs with representatives in a small neighborhood of $x_0$ in
$X.$ Denote the ideal sheaf generated by $g_1,\dots,g_n$ by $\Cal J,$ its
stalk at 0 by $J.$ Now, we use the upper--semicontinuity of the multiplicity:
$$e(\II\OO_{X(0),x_0})=e(\Cal J\OO_{X(0),x_0})\geq e(\Cal
J\OO_{X(t),x_t})\geq e(\II\OO_{X(t),x_t})=e(\II\OO_{X(0),x_0}).$$ Hence, all
inequalities are equalities and, by Rees' theorem, the ideal induced by $\Cal
J$ in $\OO_{X(t),x_t}$ is a reduction of the ideal induced by $\II$ for all
$t$ in $Y$ close to 0.

Now, note that the underlying set of the exceptional
divisor $D$ of the blowup $B$ of $X$  along $\Cal J$ equals
$Y\times\PP^{n-1}.$ Also, as the normalization $\bar B\to B$ is finite, all
components of the pull--back $\bar D$ of $D$ to $\bar B$ are equidimensional
over $Y.$ Thus, it is not hard to see, using (\class.1)(iv) that an element
induced by a section $h$ of
$\OO_X$  in
$\OO_{X,x_0}$ is integrally dependent on $\Cal J$ if, for all $t$ in a
Zariski--open subset of $Y,$ the element induced by $h$ in $\OO_{X(t),x_t}$ is
integrally dependent on the ideal induced by $\Cal J.$ In fact, the order of
vanishing of the pullback of $h$ to $\bar B$ along a component $C$ of the
underlying set of $\bar D$ can be computed in the fibre over some generic $t$
in
$Y.$ By assumption, this is no smaller than the order of vanishing of $\bar D$
along
$C.$ This proves the theorem.

We want to point out that the main part in the proof is to establish the
equidimensionality of the exceptional divisor over $Y.$

\art4 Equisingularity conditions

The Principle of specialization of integral
dependence can be used to establish criteria for equisingularity conditions.
 Since the procedure for applying the principle
is the same in all cases, we review the procedure in the case of condition $w$
for  families of hypersurfaces with isolated singularities; this case was
worked out by Teissier \cite{T1}.

Let $f:(\C^{n+2},0)\to(\C,0)$ be function, vanishing on $0\times(\C,0)$
and assume that its restrictions
$f_t=f|_{(\C^{n+1},0)\times t}$ are reduced. Let $z_0,\dots,z_n$
be coordinates on $(\C^{n+1},0)$ and $t$  the coordinate on
$(\C,0).$ Consider the hypersurface $X\subset \C^{n+2}$
defined by $f.$ By assumption, the fibres $X(t)$ over points $t$
in
$(\C,0)$ are reduced. Then, the smooth part of $X$
satisfies Verdier's condition $w$ over
$0\times(\C,0)$ at 0 if the distance between the tangent hyperplane to $X$
at a smooth point $x$ of $X$ and $0\times(\C,0)$
goes to zero no faster than the distance of $x$ to  $0\times(\C,0)$. (Here, we
use an appropriate distance function of linear subspaces of $\C^{n+1};$
see e.g. \cite {T4, Ch. 3}.) The first step is to describe this condition in
terms of integral closure.  Bernard Teissier did this by showing
that it is equivalent to the inclusion
$$\left.{\pd f\over \pd t}\right|_X\in\overline{m_zJ_z(f)\OO_{X,0}},
\tgs4.1$$ where $m_z$ is the ideal generated by the $z$--coordinates, and
$J_z(f)$ is the ideal generated by the partial derivatives of $f$ with
respect to these coordinates.

The next step is to check that the equisingularity condition holds generically.
This can be seen by using Teissier's Idealistic Bertini Theorem (see e.g.
\cite{T2, Chapter 2}). The principle of specialization then ensures that the
constancy of the multiplicity of $mJ(f_t)\sub \OO_{X(t),0}$ implies that
the equisingularity  condition holds at the origin as well.  At this point
there are two directions in which to move. One involves looking for geometric
interpretations  of the multiplicities used to control the  equisingularity
condition. In his Carg\`ese paper \cite{T1}, Teissier showed that the constancy
of
$e(mJ(f_t))$ is equivalent to the constancy of the Milnor numbers of the
plane sections of $X_t$.  The other direction is to show that the constancy of
the invariants is necessary as well as sufficient for the equisingularity
condition to hold. This can be done by showing that the equisingularity
condition gives enough control over the exceptional divisor of the  blowup of
$X$ along
$m_zJ_z(f)$ so that the constancy of the multiplicity can be deduced from the
intersection theory principle ``Conservation of Number"; see \cite{F, Ch.
10}.


\sectionhead {\segre} Segre cycles and Polar Varieties


We define the main devices of our work, Segre cycles and polar varieties of an
ideal in the local ring of a reduced analytic space germ, and discuss some
properties of their multiplicities.

Our approach was inspired by work of Teissier \cite{T4} and Massey
\cite{\Ma}. Similiar work was
done by  Leendert VanGastel \cite{Gas}, and Kleiman and Anders Thorup. In
fact, our polar
multiplicities and Segre numbers are special cases of polar multiplicities
as defined in \cite{KT1,
(8.1)}.

\stp 1 Let $(X,0)\sub (\C^N,0)$ be a reduced closed analytic space of
pure dimension $n$ and
$I\subset \OO_{X,0}$ an ideal which defines a nowhere dense subspace of
$(X,0).$ Consider the blowup
$B$ of $X$ along $I:$
$$b:B=\Bl_IX\to X$$ with exceptional divisor $D$. Suppose $I$ is generated by
$(f_1,\dots, f_M).$ Then, the blowup $B$ equals the closure in
$X\times\PP^{M-1}$ of the graph of the map
$$X-V(I)\to\PP^{M-1},x\mapsto(f_1(x):\dots:f_M(x)).$$

The restriction to $B$ of the projection onto the second factor of the
product $X\times\PP^{M-1}$ induces a map
$p:B\to\PP^{M-1}.$ Hence, a hyperplane in the projective space induces a
Cartier divisor on $B$ via the pull--back $p^*$, provided $B$ is not
contained in
the product of $X$ and the hyperplane.  We call it a {\it
hyperplane on} $B.$

Let $V$ be a reduced subspace of
$B$ of pure dimension $k$ no component of which is contained in
$D.$ A hyperplane $H$ of $B$ is {\it general with respect to} $V$
if $H\cap V$ is reduced of dimension $k-1$ and
none of its  components is contained in $D.$ Then, the intersection
equals the strict transform by $b$ of its image $b(H\cap V)$ in $X.$
In fact, outside of $D,$ the map $b$ is an isomorphism, and by assumption, the
closure of $H\cap (V-D)$ equals $H\cap V.$

Using Kleiman's Transversality Lemma, one can show that there exists a
Zariski--open subset of the set of all hyperplanes of $B$ which are
general w.r.t. $V$ (see e.g. \cite{T4,Ch.4}).

Consider a hyperplane $H$ on $B$ induced by a hyperplane $H'$ in
$\PP^{M-1}.$ An equation of
$H'$ corresponds to a linear combination
$g$ of the generators of $I.$ Suppose that $H$ is general w.r.t. $B,$ then the
above argument shows that the  topological closure of $V(g)-V(I)$ in $X$
equals the image of
$H\cap B.$ This simple observation will allows us to give an
alternative constructions of our main devices, polar
varieties and Segre cycles, in Lemma (\segre.2).

\medskip
Let $Y$ be a reduced subspace of $X$ of pure dimension $m,$  no component of
which is contained in
$V(I).$  Consider an $m$--tuple $\bold g=(g_1,\dots,g_m)$ of  linear
combinations of the generators of $I.$ Assume that each hyperplane
$H_i$ on $B$  corresponding to $g_i$ is general w.r.t.
$H_1\cap\dots\cap H_{i-1}\cap \Bl_IY.$  Then, we also say that $\bold g$ is
{\it general with respect to $Y.$} The {\it polar varieties} and {\it Segre
cycles} of
$I$ on
$Y$ with respect to
$\bold g$ are defined as follows:

$$P_0^{\bold g}(I,Y):=Y,\quad P_k^{\bold g}(I,Y):=b(H_1\cap\dots\cap
H_k\cap\Bl_IY),$$
$$ \Lambda^{\bold g}_k(I,Y):=b_*(H_1\cdots
H_{k-1}\cdot D\cdot\Bl_IY).$$

The polar varieties are reduced by the above discussion. Note, the index $k$
gives the codimension of
$P_k^{\bold g}(I,Y)$ and $\Lambda^{\bold g}_k(I,Y)$ in $Y$ if they are not
empty. For polar varieties, this follows from the definition of 'general'. A
dimension count together with the properties of the push--forward establishes
the result for Segre cycles.

\medskip
Let $I'$ be the ideal induced by $I$ in $\OO_{Y,0}$ and $\bold g'$ its set of
generators given by the restriction of the elements of $\bold g$ to $Y.$
Then, the blowup of $Y$ along $I'$  is isomorphic to the strict transform of
$Y$ by $b$, and the exceptional divisor of $Y$ is the intersection of the
blowup of $Y$ with the
exceptional divisor of $X$. Therefore the following equalities obtain for
$k=1,\dots,m-1.$
$$\align
P_k^{\bold g}(I,Y)&=P_k^{\bold g'}(I',Y),\tgs1.1 \\
\Lambda^{\bold g}_k(I,Y)&=\Lambda^{\bold g'}_k(I',Y).\tgs1.2
\endalign$$

\medskip
If $\bold g$ is formed by generic linear combinations of the generators of
$I,$ we omit the superscript $\bold g$ in the notation for polar varieties
and Segre cycles. This notation is slightly imprecise because different
generic $m$--tuples will yield, in general, different polar varieties and Segre
cycles. However, their multiplicity at 0
is independent of the choice of the generic linear combinations. Teissier
\cite{T4,IV,1.3,p.419} proved this for the special case of the Jacobian
ideal. In fact, we will see in the next section that the multiplicities of
both, generic polar varieties and generic Segre cycles,
are given by intersection formulas which are independent of the choice of
generic
$\bold g.$ Thus, the
$k$th {\it polar multiplicity} of $(I,Y)$  is well--defined for
$k=0,\dots,m-1$ as
$$m_k(I,Y)\:=\mult_0P_k(I,Y),$$ and for $k=1,\dots,m$ the $k$th
{\it Segre number} of $(I,Y)$ is $$e_k(I,Y)\:=\mult_0\Lambda_k(I,Y),$$ where
the multiplicity of a cycle
$S=\sum a_i[V_i]$ at 0 is defined as  $$\mult_0S=\sum a_i\mult_0V_i.$$
By our definition of a general, a general $m$--tuple  $\bold g$ of linear
combinations of generators of $I$ need not be generic. So, by the upper
semi--continuity of the multiplicity, we have for $k=1,dots, m-1,$
$$
m_k(I,Y)\leq\mult_0P_k^{\bold g}(I,Y)\tag \segre.1.3$$

\medskip
Let $X=(\C^{n+1},0)$ and form the Jacobian ideal  $I$ of a
function
$$f:(\C^{n+1},0)\to(\C,0)$$ generated by its partial derivatives
$$\bold g=\left({\pd f\over\pd z_0},\dots,{\pd f\over\pd z_0}\right).$$
Suppose that the coordinates on $X$ are sufficiently general so that
$\bold g$ is a general $m$--tuple of linear combinations of generators of $I.$
Then
$\Lambda_k^{\bold g}(I,X)$ is called the
$k$--codimensional L\^e cycle
$\Lambda_k^{\bold z}(f)$ of $f$ with respect to the coordinates $\bold
z=(z_0,\dots,z_n)$ (see
\cite{\Ma}). Its multiplicity
$\lambda_k^{\bold z}(f)$ at 0 is called the $k$--codimensional L\^e number of
$f$ at 0 w.r.t. $\bold z.$ We will return to this special case in Section \af.

\medskip

The following Lemma gives a useful description of Segre cycles and polar
varieties.

\lem2 Assume that $\bold g=(g_1,\dots,g_m)$ is a generic $m$--tuple of linear
combinations of generators of $I.$ For
$k=1,\dots,m-1,$ the $k$--codimensional polar variety of $I$ on $Y$ equals
the closure of
$$V(g_k|_{P^{\bold g}_{k-1}(I,Y)})-\V(I)$$
in $Y.$

Also, the following equalities of cycles obtain:

$$\eqalign{\Lambda^{\bold g}_k(I,Y)&=[V(g_k|_{P^{\bold
g}_{k-1}(I,Y)})]-[P^{\bold g}_k(I,Y)],\cr
\Lambda^{\bold g}_m(I,Y)&=[V(g_m|_{P^{\bold
g}_{m-1}(I,Y)})].}$$

\pf The first statement follows directly from the definition of the polar
varieties and the discussion that preceeded it.

For the second statement, using (\segre.1.1) and (\segre.1.2), we replace
$X$ by
$P_{k-1}(I,Y)$ and $I$ by the ideal induced by $I$ in the local ring of
$P_{k-1}(I,Y)$ at 0. Then, we may assume that $k=1.$ Note that there is
nothing to prove if $V(I)$ is of codimension bigger than one. So, we assume
that $V(I)$ has codimension one.

By definition, the one--codimensional Segre cycle of $I$ on $X$ equals the
push--forward to $X$ of the exceptional divisor of the blowup of $I$ in $X.$
Its underlying set equals the set formed by the underlying components of the
underlying set of $V(I)$ of codimension one. So, we need to show that the
multiplicity of such a one--codimensional component $W$ in the cycle
$[V(g_1)]$
equals the multiplicity of $W$ in $\Lambda_1(I,X).$

Now, the first multiplicity is, by definition, given by
the multiplicity of the ideal $(g_1)\OO_{W,X}.$ On the other hand, the
multiplicity of $W$ in $\Lambda_1(I,X)$ equals the multiplicity of the ideal
$I\OO_{W,X}$ (see \cite{\Fu, Ex. 4.3.4, p.81}). Finally, as the generators of
$I$ are assumed to be generic, the two multiplicities in question are equal
by a theorem of Samuel (see \cite{\Fu, Ex. 4.3.5(a), p.81}). The proof is now
complete.\enddemo

\medskip
The following equality follows immediately from the above
lemma:

$$e_k(I,Y)=e_1(I,P_{k-1}(I,Y)).\tgs2.1$$

Also, the lemma shows that generic polar varieties and Segre cycles do not
change if we replace
$I$ by a reduction; for instance, the ideal generated by $n$
generic linear combinations of $I.$

The lemma is a useful tool for computing Segre cycles and Segre
numbers. For example, let $I$ be a principal ideal of $\OO_{\C^n,0}$ generated
by
$f.$ Then we have $\Lambda_1(I,X)=[V(I)]$ and all other Segre cycles are
empty. Therefore, $e_1(I,X)$ equals the multiplicity of $V(f)$ at 0.

\prop3 Let $\bold g$ be an $m$--tuple of generic linear combinations
of generators of
$I,$ and
$k=1,\dots,m-1.$ Let $p:\C^N\to\C^{m-k}$ be a linear map so that the
$(N-m+k)$--plane
$H=p\inv(0)$ intersects the tangent cone at 0 of $\Lambda_k^{\bold g}(I,Y)$
transversally. Let $\epsilon$ be a general point of
$\C^{m-k}$ close to 0, and $H^\epsilon=p\inv(\epsilon).$ Then,
$$e_k^{\bold g}(I,Y)=\sum_{x\in|\Lambda_k^{\bold g}(I,Y)|\cap H^\epsilon}
e(I\OO_{P_{k-1}^{\bold g}(I,Y)\cap H^\epsilon,x}).$$

\pf By (\segre.2.1) and (\segre.1.2) we may replace $X$ by $P_{k-1}(I,Y).$
Thus, we may asssume that
$k=1,$ and
$X=Y.$ Let
$W$ be a component of the underlying set of
$\Lambda_1(I,X).$ Then, by \cite{Li, (3.5)} applied to the restriction of
$p$ to a neighbourhood in $X$ of a point $x$ in $W\cap H^\epsilon,$ we have
$$e(I\OO_{W,X})= e(I\OO_{X\cap H^\epsilon,x}).$$ By \cite{\Fu, Ex.
4.3.3, p.81}, the multiplicity of $I\OO_{W,X}$ equals the
multiplicity of $W$ in $\Lambda_1(I,X).$ Also, the number of points in
$W\cap H^\epsilon$ equals the multiplicity of $W$ at 0. This implies the
claim.\enddemo

\art4 Moving and Fixed Components

Let $Y$ be a reduced subspace of $X$ of pure dimension $m.$ A subset $W$ of
$(X,0)$ will be said to be {\it distinguished} by $(I,Y)$  if it is the image
in $Y$ of a component of $D\cdot \Bl_IY.$

A set $W_k$ of codimension $k$ in $Y$ which is distinguished by $(I,Y)$ is then
a component of $\Lambda_k(I,Y),$ as can be seen from a dimension argument.
Such a component is called a {\it fixed component} of the $k$th Segre--cycle of
$I$ on $Y.$ A component of $\Lambda_k(I,Y)$ which is not distinguished by
$(I,Y)$ is called a {\it moving component}. It is
distinguished by some $(I,P_l(I,Y))$ where $P_l(I,Y)$ is a polar variety of
dimension bigger than $\Lambda_k(I,Y).$ (It is easy to see that it is
distinguished by the pair $(I,P_{k-1}(I,Y)).$ In general it may be
distinguished by some higher--dimensional polar variety as well.)  The reason
for this terminolgy is simple. As the hyperplanes in the definition of the
L\^e cycles are varied, the fixed
components will not change, while the moving components will.

A moving component of $\Lambda_k(I,X)$ comes from an irreducible
component
$C$ of the exceptional divisor $D$ of $\Bl_IX,$ whose image  is of
codimension less than $k,$ and the fibre of $C$ over 0 has dimension at least
$k-1.$ If the dimension of this fibre equals $k-1,$ then the component $C$
will not induce moving components of the Segre cycles of $(I,X)$ of codimension
bigger than $k.$ This follows directly from the definition of Segre cycles.


\sectionhead {\form} Intersection Formulas


In (\form.2) we express the Segre numbers as intersection numbers. Hence,
intersection theory is the main tool, although only the rudiments are needed.
We review some basic definitions from Fulton's book \cite{Fu} in (\form.1). In
Lemma (\form.3) we study the relation of the blowup of an ideal $I$ on a space
germ $X$ and the blowup of the induced ideal on a hyperplane slice of $X.$
This leads to relations of the Segre numbers of $I$ and the induced ideal on
the hyperplane slice, as described in (\form.4). The expansion formulas
(\form.5) express the Segre numbers of the product of $I$ and the maximal
ideal in terms of the Segre number of $I$ and its polar multiplicities.
We restate a result of Kleiman and Thorup that gives a length theoretic
interpretation of the top Segre number in (\form.7).

\art 1 Intersection Theory

Let $X\sub (\C^{N},0)$ be a reduced analytic
space germ of pure dimension $n$   containing 0, let $I\sub\OO_{X,0}$ be an
ideal and $m\subset\OO_{X,0}$ be the maximal ideal.  Let
$\Lambda_k(m\,I,X)$ be a generic $k$--codimensional Segre cycle of
$m\,I$ and $e_k(m\,I,X)$ its multiplicity at 0. The main tool for studying
these Segre cycles is the following commutative diagram. $$\CD  D\subset@.
B=\Bl_{m\,I}X@>\tilde b_1>> B_2=\Bl_{I}X@.\supset D_2\\ @.@V\tilde b_2VV @V
b_2 VV@.\\ D_1\subset@. B_1=\Bl_mX@>b_1>>X@.\\ \endCD$$

The exceptional divisors of the blowups $B_1,B_2$ are $D_1,D_2.$ The
exceptional divisor of the blowup $b:\Bl_{m\,I}X\to X$ is denoted $D.$

We will use the first Chern classes of the tautological line bundles on the
blowups  $h_1=c_1(\OO_{B_1}(1)),h_2=c_1(\OO_{B_2}(1)),h=c_1(\OO_{B}(1)).$ The
reader should think of such a class, say $h_1$, as an operator that gives,
applied to an irreducible variety $V\sub B_1\subset (X\times\PP^n)$, the
rational equivalence class of $V\cap H,$ where $H$ is a {\it generic
hyperplane} of $B_1$, that is, $H=(X\times H')\cap B_1$ with $H'$ a generic
hyperplane of $\PP^n.$  The divisor $H$ represents the tautological line
bundle on $B_1.$  (The correctness of this view point follows from Kleiman's
Transversality Lemma which implies that $H$ intersects $V$ transversally
and the intersection is reduced and of dimension $\dim V-1$.) This
definition extends to cycles by linearity.

The line bundle on $B_i$
associated to $D_i$ is the dual of $\OO_{B_i}(1).$ Hence, intersecting with
$D_i$ equals the operation of $-h_i.$ We say that (the line bundle
associated to) $D_i$ is the dual of $\OO_{B_i}(1).$

For $i=1,2$ the pull--back of $h_i$ to $B$ will be denoted by $\th1$ and the
pull--back of $D_i$ in $B$ by $\tD{i}.$ Then, the following equalities of
Cartier divisors, resp. first Chern classes, obtain (see \cite {KT2, (2.7)})
$$D=\tD1+\tD2,\quad h=\th1+\th2.\tag
\form.1.1$$

For a cycle $S$ on $B$ the part of $S$ formed by the components whose
generic points map into $X-0$ will be denoted $S^{X-0}.$ In other words, the
cycle
$S^{X-0}$ is formed by those components not lying in the fibre of $B$ over 0.
Note that for a generic hyperplane $H$ on $B$ we have  $$H\cap S^{X-0}=(H\cap
S)^{X-0}.\tag\form.1.2$$ To see this we count dimensions, assuming $H\cap
S^{X-0}$ is not empty. Then, by Kleiman's Transversality Lemma, we may assume
that $H$ intersects every component of the fibre of $S^{X-0}$ over 0 properly.
Hence,
$$\dim H\cap S^{X-0}=\dim  S^{X-0}-1\geq
\dim( S^{X-0}\cap b\inv(0)) > \dim (H\cap S^{X-0}\cap b\inv(0)).$$ Therefore, no
component of the fibre of
$H\cap S^{X-0}$ over 0 is a component of $H\cap S^{X-0}$.

Clearly, the same argument works for hyperplanes representing $\th1$ and
$\th2.$

Also, denote the part of $S$ formed by the components mapping into 0 by $S^0.$
Then, (\form.1.2) is equivalent to $$H\cap S^0=(H\cap S)^0.\tag\form.1.3$$

We will identify a Cartier divisor $D$ with its associated Weil divisor. We
also write
$D^0$ for the part of its associated Weil divisor formed by the components
mapping to 0. The cycle
$D^{X-0}$ is defined analogously.

\medskip
The {\it degree} of a cycle $S$ is denoted $\int S.$ It is the sum
of the multiplicities of the 0--dimensional components. For a
cycle $S$ in the fibre of $B$ over 0, the degree of $S$ depends only on the
rational equivalence class of $S$ in this fibre. For cycles outside this
fibre the degree is no invariant of rational equivalence.

\art 2 Intersection Formulas

We can express the multiplicities of the Segre cycles at 0 as intersection
numbers   $$\align
e_k(m\,I,X)&=\int \th1^{n-k-1}h^{k-1}D^{X-0}\cdot\tD1,
\quad k=1,\dots,n-1,\tag\form.2.1\\
e_n(m\,I,X)&=\int h^{n-1}D^0,\tag
\form.2.2\\
e_k(I,X)&=\int\th1^{n-k-1}\th2^{k-1}\tD2^{X-0}\cdot\tD1, \quad
k=1,\dots,n-1,\tag \form.2.3\\
e_n(I,X)&=\int h_2^{n-1}D_2^0.\tag\form.2.4\\
\endalign$$

The following formula is a generalization of L\^e and Teissier's polar
multiplicity formula \cite{LT2, (5.1.1)}
$$m_k(I,X)=\int\th1^{n-k-1}\th2^k\tD1.\tgs2.5$$
They considered the special case of the Jacobian ideal, but their proof works
also in our setup. Alternatively, the formula can be verified by using
similar arguments to the ones that will imply (\form.2.1)--(\form.2.4); see
also \cite{KT1, (8.9)}.

We extend the notion of strict transform to cycles by linearity. Then the
blowup of a cycle is defined to be its strict transform
by the blowup--map.

To prove the first formula consider
$$\Bl_m\Lambda_k(m\,I,X)=\Bl_m(b_*(D\cap H_1\cap\dots\cap H_{k-1}))=
\tb{2*}(D^{X-0}\cap H_1\cap\dots\cap H_{k-1}),$$ where the $H_i$ are generic
hyperplanes of $B.$ The first equality
follows from the defintion of
$\Lambda_k(m\,I).$ By the argument that led to (\form.1.2) only the components
of $D$ mapping to a subset of codimension at least $k$ are not annihilated by
intersecting with $k-1$ generic hyperplanes of $B.$ The second equality
follows.  Now,
$$\eqalign{e_k&(m\,I,X)=\int h_1^{n-k-1}\Bl_m\Lambda_k\cdot D_1\cr &=\int
\th1^{n-k-1}(D^{X-0}\cap H_1\cap\dots\cap H_{k-1})\cdot\tD1
=\int\th1^{n-k-1} h^{k-1}D^{X-0}\cdot\tD1.}$$ Here, the first equality is
the well--known intersection formula for the multiplicity of a cycle. (In
fact, it is just a translation of Samuel's definition of the multiplicity of
an ideal into geometry.) The second equality follows from the equality above
and  the projection formula \cite{F, Prop. 2.5, p.41}. Then, we pass to
rational equivalence to get the formula.

The second formula follows directly from the definition of $e_n$ and
(\form.1.3).  The proof of the next two formulas runs analogously.

\lem 3 For $k=1,\dots,n-1$ let $L_{n-k}$ be a generic $(n-k)$--codimensional
linear subspace of $\C^N.$ Consider the blowup $\Bl_I(X\cap L_{n-k})$ with
exceptional divisor $D_{2,n-k}.$ Then, its pull--back $\tD{2,n-k}$ to $B$
satisfies the following relations up to rational equivalence:
$$\align
\tD{2,n-k}^0&=\th1^{n-k}\tD2^0=\th1^{n-k-1}\th2\tD1+
\th1^{n-k-1}\tD2^{X-0}\cdot\tD1,\\
\tD{2,n-k}^{X-0}&=\th1^{n-k}\tD2^{X-0}.\endalign$$

\pf
We will only show the formulas for $L_1.$ The general case follows by
induction. Now,
$$\Bl_{mI}(X\cap L_1)=\Bl_{b_1^*I}(H\cap B_1)=(\tb2\inv H)\cap B,$$ where $H$
is the hyperplane on $B_1$ induced by the hyperplane in $\PP^{N-1}$
corresponding to $L_1.$ Indeed, the first equality follows from (\form.1.3),
and the second from the properties of the pull--back of a line bundle.
Passing to rational equivalence implies the second relation. For the first
we consider
$$\tD{2,1}^0=\th1\tD2^0=(-\tD1)\cdot(\tD2-\tD2^{X-0})=\th2\tD1+
\tD2^{X-0}\cdot\tD1.$$ Here, the first equality follows from the
properties of $\th1,$ the second and third use duality. \enddemo

\art 4 Hyperplane Sections and Segre Numbers

Let $L_{n-k}$ be a generic $(n-k)$--codimensional linear subspace of
$\C^N.$ Then there are the following relations of Segre numbers of $(I,X)$
and those of $(I,X\cap L_{n-k}).$

$$\align
e_i(I,X\cap L_{n-k})&=e_i(I,X)\quad\text{for }i=1,\dots,k-1,\tag\form.4.1\\
e_k(I,X\cap L_{n-k})&=m_k(I,X)+e_k(I,X).\tag\form.4.2
\endalign$$
Indeed, the same argument as above shows that the divisor $\tD{1,n-k}$,
defined by the pull--back of the maximal ideal $m$ to $\Bl_{mI}(X\cap
L_{n-k})$ is given by $\th1^{n-k}\tD1.$ Then the first relation follows
immediately from the above lemma (\form.3) and the intersection formula
(\form.2.3):
$$\eqalign{e_i(I,X\cap
L_{n-k})&=\int\th1^{k-i-1}\th2^{i-1}\tD{2,n-k}^{X-0}\tD{1,n-k}\cr
&=\int\th1^{n-i-1}\th2^{i-1}\tD2^{X-0}\tD1=e_i(I,X).}$$

For the second relation observe that the push--forward of $\tD{2,n-k}^0$
to $B_2$ equals $D_{2,n-k}^0.$ Hence, we can use the projection formula
and (\form.3) to compute $e_k(I,X\cap L_{n-k})$ on $B:$   $$e_k(I,X\cap
L_{n-k})=\int\th2^{k-1}(\th1^{n-k-1}\th2\tD1+
\th1^{n-k-1}\tD2^{X-0}\cdot\tD1).$$ The polar multiplicity
formula (\form.2.5) and (\form.2.3) yield the desired relation.

\thmx{The expansion formulas}5
In the setup (\form.1) the following formulas obtain.
$$\align e_n(m\,I,X)&=\sum_{i=0}^{n-1} \binom n i m_i(I,X)+\sum_{i=1}^n
\binom{n-1}{i-1}e_i(I,X) ,\tag\form.5.1\\
e_k(m\,I,X)&=\sum_{i=1}^k
\binom{k-1}{i-1}e_i(I,X).\tag\form.5.2 \endalign$$

\pf Using the above formula (\form.1.1), we expand $h=\th1+\th2$ in
(\form.2.1), and get
$$e_k(m\,I,X)=\sum_{i=0}^{k-1}\binom{k-1}i\int
\th1^{n-i-2}\th2^i\tD2^{X-0}\cdot\tD1=\sum_{i=0}^{k-1}\binom{k-1}i
e_{i+1}(I,X).$$   This implies the second formula.

For the first formula we also expand $h=\th1+\th2$ and $D^0=\tD1+\tD2^0,$
and use the Polar Multiplicity formula (\form.2.5):
$$\eqalign{e_n(m\,I,X)&=
\sum_{i=0}^{n-1}\binom{n-1}i\int\th1^{n-i-1}\th2^i(\tD1+\tD2^0)\cr
&=\sum_{i=0}^{n-1}\binom{n-1}i (m_i(I,X)+\int\th1^{n-i-1}\th2^i\tD2^0).}$$
Now, by (\form.3) and (\form.4.2), for $i$ less than $n-1$ the degree of
$\th1^{n-i-1}\th2^i\tD2^0$ equals $e_{i+1}(I,X)+m_i(I,X).$ Hence,
$$\eqalign{e_n(m\,I,X)=\sum_{i=0}^{n-1}&\binom{n-1}i m_i(I,X)\cr
&+\sum_{i=1}^{n-1}\binom{n-1}{i-1}
\left(e_{i}(I,X)+m_{i}(I,X)\right)+e_n(I,X).}$$ Adding the two sums yields
the formula. \enddemo

\rmk 6  The results of this section remain valid, if the maximal ideal $m$
is replaced by an $m$--primary ideal $m'.$ This should be useful in the study
of weighted homogeneous polynomials: Let $m'=(z_1^{p_1},\dots,z_N^{p_N})$
where the $p_i$'s are positive integers. We define the weight of $z_i$ to
be $1/p_i.$ Let
$X=V(f)$ be a hypersurface given by a weighted homogeneous polynomial $f$ of
weighted degree $d.$ Then we have
$m_0(m',X)=d\,p_1\dots p_N$ (see \cite{\Li, p.113} and \cite{\MO, p.387}).

In fact, it is possible to replace $m$ by an ideal $m'$ that induces an
ideal of finite codimension in the local ring of $V(I)$ at 0 but not on
$\OO_{X,0}.$ Still, most of the results of this section remain valid. This case
seems to play an important role in the study of Thom's $a_f$ condition on
singular spaces.

\medskip
In \cite{KT2,(3.6)} Kleiman and Thorup define the generalized Samuel
multiplicity of a closed subscheme on a module at a closed point. We apply
their proposition (3.5) to the subscheme $Z$ defined by $I$ in $X,$ the
structure sheaf of $X$ and the origin. The $k$th infinitesimal neighborhood of
0 in $Z$ is denoted $Z_k.$

\prop 7 As a function of $m,$ the dimension
$$p(m)=\dim(\oplus_{i=0}^m(I^i/I^{i+1})\otimes\OO_{Z_k})$$ is
eventually a polynomial of degree at most $r=\dim Z.$ Moreover, for $k\gg 0,$
the coefficient of $m^r/r!$ is independent of $k$ and equals $e_n(I,X).$

\pf Use \cite{KT2,(4.3)} to see the connection of the statement of
\cite{KT2,(3.5)} with the Segre numbers of the ideal $I.$
\enddemo

This shows that the top Segre number of an ideal generalizes in a natural way
Samuel's multiplicity. No such length--theoretic interpretation is known for
lower dimensional Segre numbers.


\sectionhead {\sid} The Specialization of Integral Dependence

This section contains our main result, the principle of specialization of
integral dependence. We consider a family of analytic spaces and a given ideal
sheaf on the total space. As a start we study the first Segre numbers of the
ideals induced on the fibers by the given one on the total space. They are
upper semi--continuous. The first Segre number controls the components of
the divisor mapping to codimension one subsets of the special member of the
family; see (\sid.4). Proposition (\sid.3) shows that if the first Segre
numbers are constant then the special fibre of the one--codimensional polar
variety of the total space equals the one--codimensional polar variety of the
special fibre. This serves  as the starting point of an induction in the proof
of the lexicographically upper semi--continuity of the Segre numbers
(\sid.5). The next proposition (\sid.6) links the constancy of the Segre
numbers in a family to a very strong equidimensionality condition. Finally,
the principle of specialization of integral dependence (\sid.7) follows from a
classical connection of integral dependence with this exceptional divisor.

As a corollary of the principle of specialization we
obtain a generalization of Rees' theorem (\class.2). For an ideal whose
co--support is nowhere dense we show that its Segre numbers control its
integral closure. Remark (\sid.11) indicates that a similar result is true for
arbitrary ideals when the Segre numbers are replaced by numbers that come from
the completed normal cone of the ideal.

As a first application we give a numerical characterization of limiting
tangent hyperplanes of a hypersurface.

\stp 1 Let $F:(X,0)\to(Y,0)$ be a map of germs of analytic spaces. Assume that
the fibres $X(t)$ are reduced and equidimensional of the same dimension $n$
at least 1. We assume that $X$ is embedded in $(\C^N,0)\times(Y,0)$ and that
$F$ is induced by the projection onto the second factor. Furthermore, assume
that
$0\times(Y,0)$ is contained in $X;$ we will identify it with $Y.$ Its
defining ideal sheaf in $X$ will be denoted $m_Y.$

Let $I\sub\OO_X$ be a sheaf of ideals on $X$ and denote the ideal induced
by $I$ in $\OO_{X(t),0}$ by $I(t).$ Assume that $Z(t)=V(I(t))$ is nowhere
dense in $X(t).$

We work in a modified version of the setup (\form.1) where $m$ is replaced by
$m_Y.$ In addition, for $t$ in $Y,$ we consider the blowup
$B_{2,t}=\Bl_{I(t)}X(t)$ with exceptional divisor $D_{2,t}.$ It sits inside
the fibre $B_2(t)$ of $B_2$ over $t,$ and the tautological bundle
$\OO_{B_2}(1)$ restricts to the tautological bundle $\OO_{B_{2,t}}(1).$ We
define $B_t$ and $B_{1,t}$ in the analogous way.

The Segre numbers $e_k(I(t),X(t))$ will be denoted by $e_k(t),$ for
$k=1,\dots,n.$

\prop 2 The map $t\mapsto e_1(t)$ is upper
semi--continuous.

\pf (i) {\it The map
$t\mapsto e_1(t)$ is `non--decreasing'} ; that is, if
$A$ is an analytic subvariety of $Y$ and $t$ is a general point of $A$ close
to 0, then $e_1(t)\leq e_1(0).$

As our numbers are defined on the fibres, we may replace $A$ by $Y.$ We are
going to use the classical result on the upper semi--continuity of
multiplicities of ideals of finite colength; see e.g.
\cite{Li,(3.5)}. Choose a linear map $p:\C^N\to\C^{n-1}$ the kernel of which
intersects the tangent cone of $\Lambda_1(I(0),X(0))$ at 0 transversally. Next,
we may choose a neighborhood $U$ of 0 in $\C^N$ and a representative of $(Y,0)$
which we denote again by $Y$ so that, after identifying all germs involved
with their representatives on $U\times Y$, the map
$$\pi:\Lambda_1(I,X)\to(\C^{n-1}\times Y),(z,t)\mapsto (p(z),t)$$ is finite.
Also, we may assume that the map
$$(\epsilon,t)\mapsto \sum_{x\in p\inv(\epsilon)\cap
X(t)}e(I\OO_{p\inv(\epsilon)\cap X(t),x})\tgs 2.1$$ is upper
semi--continuous on the image of $\pi.$

Next, choose general points $t$ in $Y$ and $\epsilon$ in $\C^{n+1}$ close to 0
so that
$$e_1(0)=\sum_{x\in p\inv(\epsilon)\cap
X(0)}e(I\OO_{p\inv(\epsilon)\cap X(0),x})$$ and so that the above map
(\sid.2.1) attains its generic value at $(\epsilon,t).$ Finally, again by upper
semi--continuity, we may assume that this
number is at least
$e_1(t).$ (Here, we need to consider the family of finite
projections of $\Lambda_1(I(t),X(t))$ onto $\C^{n-1}$.) This proves the
`non--decreasing' statement.

(ii) {\it The map $t\mapsto e_1(t)$ is constant on an Zariski--open subset of
$Y.$} Consider a component $C$ of $\tD2^{X-Y}$ the fibre of which over 0 in
$Y$ maps to a codimension one  subset in $X(0).$ In general, not all
components of the intersection
$C\cap \tD1$ will map onto $Y.$ For a point $t$ in the image of a
component of this intersection that doesn't map onto $Y$ we don't expect
$e_1(t)$ to be generic. Therefore, consider the Zariski--closed subset
$F$ in
$Y$ formed by the images of such `vertical' components, and, in addition,
by the singular locus of the underlying set of
$Y$ and the images in
$Y$ of components of $\tD1$ that don't map
onto
$Y.$ Then, for $t$ in $Y-F$, a dimension argument shows that the intersection
with
$\tD{1,t}$ of the part of $\tD{2,t}$ formed by components that map to subsets
of codimension one in $X(t)$ equals the fibre over $t$ of the part of
$\tD2\cdot\tD1$ formed in the same way; see the proof of \cite{GK, (2.1)}.
This is the part of the intersection that is relevant to the computation of
$e_1(t).$ It follows that
$$e_1(t)=\int\th1^{n-2}\tD{2,t}^{X(t)-0}\tD{1,t}=
\int\th1^{n-2}(\tD2^{X-Y}\tD1)(t)$$ is independent of $t$ in $U$ by
`conservation of numbers'; see \cite{F, Prop.10.2, p.180} and the discussion
in the proof of (\sid.6).
\enddemo

\prop3 Assume that the map $t\mapsto e_1(t)$ is constant. Then, the
one--codimensional polar variety of $I$ on $X$ specializes. That is,
$$P_1(I,X)(0)=P_1(I(0),X(0)).$$

\pf By Lemma (\segre.2), the polar variety $P_1(I(0),X(0))$ of the fibre
$X(0)$ is contained in the fibre of $P_1(I,X)$ over 0. So, we have to show
that no component of the fibre of
$P_1(I,X)$ over 0 is contained in the fibre of $V(I)$ over 0.

To see this, we choose, as in the proof of (\sid.2), a generic linear
projection $p:\C^N\to\C^{n-1}$ Then, by a similar argument as in the above
proof and the assumptions, the map (\sid.2.1) has constant value $e_1(0).$
Hence, for a general
$\epsilon$ in
$\C^{n-1}$ close to 0, the multiplicity of the ideal $I'$ induced by $I$ on
the fibres of the family of curves $X'=(p\inv(\epsilon)\times Y)\cap X$ is
constant. Also, note that $P_1(I',X')=P_1(I,X)\cap(p\inv(\epsilon)\times Y).$
Hence, it is enough to show that $P_1(I',X')$ is empty. This follows from
the proof of the classical Principle of Specialization; see (\class.3). In
fact, the constancy of the multiplicity of $I'\OO_{X'(t)}$ implies that each
component of the exceptional divisor of the blowup of $X'$ along $I'$ is
mapped onto $Y$ by $F\circ b_2.$ The desired result follows.
\enddemo

\rmk4 The above proof shows that the map $t\mapsto e_1(t)$ controls components
of the exceptional divisor $D_2$ the fibres of which over $0$ in $Y$ map to
codimension one subsets in $X(0).$ If it is constant then all such components
 are mapped onto $Y$ by $F\circ b_2.$

Also, the constancy of $e_1(t)$ implies that the image of such a  component in
$X$ contains $Y.$ If not, the generic value of the map (\sid.2.1) would be
bigger than $e_1(0).$

\cor5 The $n$--tuple $(e_1(t),\dots,e_n(t))$ is lexicographically  upper
semi--continuous: If for some $k$ the map $t\mapsto(e_1(t),\dots,e_k(t))$ is
constant on $Y,$ then $t\mapsto e_{k+1}(t)$ is upper semi--continuous. In
particular, $e_1(t)$ is always upper semi--continuous.

\pf  The existence of a Zariski--open subset of $Y$ on which $t\mapsto
e_{k+1}(t)$ is constant follows, as in the proof of (\sid.2), from
conservation of numbers. (See also the first part of the proof of (\sid.6).)

We are going to prove by induction
on $k$ that the constancy of the map $t\mapsto (e_1(t),\dots,e_k(t))$ implies
$e_{k+1}(0)\geq e_{k+1}(t)$ for all $t,$ and the
equality $P_k(I,X)(0)=P_k(I(0),X(0)).$  The case $k=0$ is proven in (\sid.2)
and (\sid.3). So, assume that $k$ is non--zero.

Fix a generic $n$--tuple $\bold g$ of linear combinations of generators of
$I.$ Then, by the induction hypothesis, we have
$$e_k(0)=e_1(I(0),P^{\bold g}_{k-1}(I,X)(0)),\quad\text{and}\quad e_k(t)\leq
e_1(I(t),P^{\bold g}_{k-1}(I,X)(t))
$$ for any $t$ in $(Y,0)$ by (\segre.1.4) and (\form.2.1). Now, by assumption,
the map
$t\mapsto e_k(t)$ is constant. Hence, the upper semi--continuous map
$$t\mapsto e_1(I(t),P^{\bold g}_{k-1}(I,X)(t))$$ is also constant. Therefore,
by (\sid.3), $P_k^{\bold g}(I,X)=P_1(I,P_{k-1}^{\bold g}(I,X))$
specializes. Furthermore, by (\sid.2), we have for $t$ in $Y$ close to 0
$$e_{k+1}(t)\leq e_1(I(t),P^{\bold g}_k(I,X)(t))\leq e_1(I(t),P^{\bold
g}_k(I,X)(0))=e_{k+1}(0).$$ This finishes the proof of the corollary.
\enddemo

\prop6 (i) If the map $t\mapsto (e_1(t),\dots,e_n(t))$ is constant on $Y,$ then
for a component $C$ of $D_2$ all components of
$C(Y)=C\cap b_2\inv(Y)$ are equidimensional over $Y.$

(ii) If, in addition, 0 is a regular point of $Y$ then the constancy of the
map in (i) is equivalent to all components of $D_2^Y$ and $\tD2^{X-Y}\cdot\tD1$
being equidimensional over
$Y.$

\pf (i) We do induction on $n,$ the
dimension of the fibres of $X$ over $Y.$ The case
$n=1$ is the classical case of an ideal of finite co--support over $Y.$

Now, assume that the assertion holds for families of $(n-1)$--dimensional
spaces. Let $X\to Y$ be a family of fibre dimension $n$ as in the above
setup. If for a component $C$ of $D_2$ the intersection $C\cap b_2\inv(0)$
is zero--dimensional, then, by Remark (\sid.4), its image in $X$
contains $Y.$ Hence, the map $C(Y)\to Y$ is equidimensional. For other
components $C$ of $D_2$ the claim holds if, and only if, it holds for
the intersection of $C$ with a generic hyperplane of $B_2$; see
(\form.1.2). Therefore, by the induction hypothesis, it is enough to show
that the numbers
$$e_1(I(t),P_1(I,X)(t)),\dots,e_{n-1}(I(t),P_1(I,X)(t))$$
are independent of $t$ in $Y.$

We will show this by induction on the codimension of the Segre numbers. So,
assume we have shown the claim for the numbers of codimension $1$ to $k-1.$
Then, we have
$$\CD
e_{k+1}(I(t),X(t))@.=e_{k+1}(I(0),X(0))@.=e_k(I(0),P_1(I(0),X(0)))\\
||@.@.||\\
e_k(I(t),P_1(I(t),X(t)))@.\leq e_k(I(t),P_1(I,X)(t))@.\leq
e_k(I(0),P_1(I,X)(0))\\ \endCD$$
The vertical equality on the left hand side follows from (\segre.2.1), the
other vertical equality follows from (\sid.3). In the top row, the second
equality is (\segre.2.1), while in the bottom row the first inequality
follows from (\segre.1.4), and the second from the upper semi--continuity
(\sid.5). It follows that
$e_k(I(t),P_1(I,X)(t))=e_{k+1}(t)=e_{k+1}(0)$ for all $t.$ This finishes the
proof of claim (i).

(ii) Assume that the map in (i) is constant, and consider a
component $\tilde C$ of $\tD2^{X-Y}$ and its image $C$ in $D_2^{X-Y}.$ Then,
by (i), the map $C(Y)\to Y$ is equidimensional, say of fibre dimension $k.$
Let $C'$ be the intersection of $C$ with $k$ generic hyperplanes
of $B_2.$ Then the map $C'(Y)\to Y$ is finite, and the assumption implies
that the multiplicity of $(b_{2*}C')(t)$ at 0 is independent of $t.$ As $Y$ is
smooth at 0, the ideal induced by $m_Y$ in $\OO_{X(t),0}$ equals the maximal
ideal of this local ring. Hence, by a projection formula for multiplicities
\cite{F,(4.3.6)}, this multiplicity equals the multiplicity of the ideal
induced by $b_2^*(m_Y)$ on
$C'(t).$ Now, by the classical Principle of Specialization, the
equimultiplicity implies that the exceptional divisor of the blowup of $m_Y$
in $C'$ is equidimensional over
$Y.$ Therefore, the intersection $\tilde C\cap\tD1$ is equidimensional over
$Y;$ and, therefore, so is $\tilde C(Y)\to Y.$

Now, assume that the equidimensionality condition holds, and that the
underlying set of $Y$ is smooth at 0. Then the dimension argument
\cite{GK, (2.1)}, together with the commutativity of refined Gysin
homomorphisms, shows:
$$D_2^Y(t)=D_{2,t}^0,\qquad
(\tD2^{X-Y}\cdot\tD1)(t)=\tD{2,t}^{X(t)-0}\cdot\tD{1,t}.$$  (Let $t\buildrel
i\over
\hookrightarrow Y$ be the regular embedding. The
equidimensionality assumptions and the
compatibility properties \cite{F, 6.2.1} of the refined Gysin
homomorphism $i^!$ imply
$i^!(\tD2^{X-Y}\cdot\tD1)=(\tD2^{X-Y}\cdot\tD1)(t)$ and
$\tD{2,t}^{X(t)-0}\cdot\tD{1,t}=\tD{2,t}^{X(t)-0}\cdot\tD1.$) Therefore, we
have for
$k=1,\dots,n-1$ that
$$e_k(t)=\int\th1^{n-k-1}\th2^{k-1}\tD{2,t}^{X(t)-0}\tD{1,t}=
\int\th1^{n-k-1}\th2^{k-1}(\tD2^{X-Y}\tD1)(t)$$ is independent of $t$ by
`conservation of numbers'. (The smoothness of $|Y|$ is needed to apply
`conservation of numbers'.) A similiar argument shows that the map
$t\mapsto e_n(t)$ is constant.
\enddemo

\medskip The equidimensionality of $\tD2^{X-Y}\cdot\tD1$ over $Y$ implies the
equidimensionality of $D_2^{X-Y}\cap b_2\inv(Y).$ The converse doesn't hold
in general: Consider a family of reduced plane curves over $(\C,0)$
given by a function $f:(\C^3,0)\to(\C,0).$ Let $I$ be the ideal generated by
$f.$ The blowup of the principal ideal $I$ is
isomorphic to $(\C^3,0),$ the divisor $D_2$ is given by the vanishing of $f.$
By construction, the intersection of $D_2^{X-Y}=D_2=V(f)$ with the parameter
axis is trivially equidimensional over $(\C,0).$ But, the one--codimensional
Segre number
$e_1(t)$  equals the multiplicity at 0 of the restriction of $f$ to $t\times
(\C^2,0)$ which need not be independent of $t.$

\thmx{The Principle of Specialization}7 In the setup (\sid.1) let
$h\in\OO_{X,0}$ be a function so that $h\OO_{X(t),0}$ is integrally dependent
on $I(t)$ for all $t$ in a Zariski--open subset of $Y.$ If the numbers
$e_1(t),\dots,e_n(t)$ are independent of $t,$ then $h$ is integrally
dependent on $I.$

\pf The proof proceeds similarly to Teissier's original proof \cite{T3,
Appendice I}. By the above Proposition (\sid.5) the assumption implies that all
components of
$D_2$ map onto
$Y.$ Therefore, the same is true for the exceptional divisor $\bar D_2$ of the
normalized blowup of $X$ along $I.$ Now, we use the characterization
(\intro.1)(iv) of integral dependence. The fibre $C(Y)$ of a component $C$
of
$\bar D_2$ over $Y$ maps onto $Y.$ Hence, we can compute the order of
vanishing of the pullbacks of $h$ and $I$ along $C$ at points that map to the
Zariski open subset of $Y$ in the statement of the theorem. Then, it is not
hard to see that, after shrinking this open subset $U$ of $Y,$ for a
point $t$ in $U$ the order of vanishing along $C(t)$ of the
pullback of $I(t)$ in the normalized blowup of $X(t)$ along $I(t)$ equals
the order of vanishing of the pullback of $I$ along $C;$ and similarly for
$h$ (see e.g.
\cite{T4, Ch. I, 1.3.4 and 1.3.6}). This proves the claim.\enddemo

\rmk8 It is easy to see that a similar statement obtains if the embedding of
$Y$ into $X$ is replaced by a subspace $S$ in $X$ that is finite over $Y$:
We replace the numbers $e_k(t)$ by the sum of the Segre numbers of
the ideal induced by $I$ in the local rings $\OO_{X(t),x^i_t}$ where the
points $\{x^i_t\}$ form the underlying set of $S(t).$ Note, however, that the
analog of (\sid.6)(ii) is false. In fact, the defining ideal of $S$ may not
induce the maximal ideal in $\OO_{X(0),0}.$

\medskip
As a corollary of the principle of specialization of integral dependence we
obtain the following generalization of Rees' Theorem (\intro.2).

\cor 9 Let $X\sub(\C^N,0)$ be an analytic germ of pure dimension $n$ and
$I\subset J\subset\OO_{X,0}.$ Then, $\bar I=\bar J$ if, and only if,
$e_k(I)=e_k(J)$ for $k=1,\dots,n.$

\pf Consider the family $X'=X\times(\C,0)\to(\C,0)$ and the ideal
$I'=(I+tJ)\OO_{X'},$ where $t$ is the coordinate on $(\C,0).$
Then, $I'(0)=I$ and $I'(t)=J$ for $t\neq 0.$ Now, if their Segre numbers
are equal, the theorem  implies $\overline{I\OO_{X'}}=
\overline{J\OO_{X'}},$ and so $\bar I=\bar J.$

The other direction follows from the existence of a generically
one--to--one map $\Bl_JX\to\Bl_IX,$ which is compatible with the exceptional
divisors of the blowups, if
$\bar I=\bar J.$ \enddemo

\rmk {10} B\"oger's theorem (\cite{B}, \cite{Li}) is an easy corollary of
this last result.
The setting of B\"oger's theorem is an ideal $I$ in an equidimensional
local ring $R$ with the
property that $I$ has a reduction generated by the same number of
generators as the height of $I$.
This implies that $V(I)$ is equidimensional. B\"oger's theorem says that if
$J$ is an ideal
with the same radical as $I$, $I\subset J$ then $J$ is in the integral
closure of $I$ iff the
multiplicity of
$J$ in each of the local rings $R_P$ is the same as the multiplicity of $I$
in $R_P$, $P$ varying
through the minimal primes of $I$. The hypothesis on $I$ implies that the
only non-zero Segre number
is
$e_j(I)$, where
$j$ is the height of $I$; the hypothesis on $J$ implies that
$e_i(J)=e_i(I), i\le j$. The
lexicographic upper semicontinuity of the
the Segre numbers then implies that $e_i(J)=e_i(I)=0, i>j$. Then our
extension of Rees' theorem
implies that ${\overline I}={\overline J}$.

\art11 Ideals with dense co--support

Let $(X,0)$ be a reduced analytic space germ, and $I\subset\OO_{X,0}$ an
ideal. We don't require that $V(I)$ is nowhere dense in $X.$ Following ideas
of Kleiman and Thorup's work \cite{KT2}, we consider the embedding and
projection
$$(X,0)=(X,0)\times 0\hookrightarrow(X\times\C,0)
\buildrel p\over\to (X,0).$$ Let
$y$ be the coordinate function on $(\C,0).$ To an ideal $I$ of $\OO_{X,0}$ we
associate the ideal $\hat I$ of $\OO_{X\times\C,0}$ generated by the
pull--back
$p^*I$ and $y$ in $\OO_{X\times \C,0}.$ Let $h$ be an element of $\OO_{X,0}.$
Then,
$h$ is integrally dependent on $I$ if, and only if, its pull--back $p^*h$ is
integrally dependent on $\hat I.$ In fact, an integral relation
$$h^k+a_1h^{k-1}+\dots+a_k=0,\quad a_j\in I^j,$$
pulls back to an integral relation in $\OO_{X\times\C,0}$ with coefficients
in the correct powers of $\hat I.$ On the other hand, an integral relation
$$(p^*h)^k+a'_1(p^*h)^{k-1}+\dots+a'_k=0,\quad a'_j\in \hat I^j,$$
restricts to the required relation for $h$ in $\OO_{X,0}$ with coefficients
in the correct powers of $I.$

It follows that Theorem (\sid.3) holds in this more general setup if the
Segre numbers of $I$ on $(X,0)$ are replaced by the Segre numbers of $\hat I$
on
$(X\times \C,0).$ In fact, we can compute the Segre numbers of $\hat I$ by
using the pull--back $p^*(m)$ of the maximal ideal of $\OO_{X,0};$ see also
Remark (\form.6).

\art 12 Limiting Tangent Hyperplanes

Consider a reduced hypersurface $X=V(f)\subset(\C^{n+1},0).$ A hyperplane $H$ is
a {\it limiting tangent hyperplane} of $X$ at 0 if it is the limit of tangent
hyperplanes of $X$ at smooth points converging to 0. There is a criteria for
limiting tangent hyperplanes in terms of  integral closure. Let
$J(f)\subset\OO_X$ be the Jacobian ideal. It is generated by the partial
derivatives of $f.$ Consider the ideal $J(f)_H\sub J(f)$ generated by the
directional derivatives $\pd f\over\pd v$ with $v$ in $H.$ Then, $H$ is
a limiting tangent hyperplane if, and only if, $J(f)_H$ is not a reduction of
$J(f)$ (see \cite{T1, p.321,308}). Using the above Corollary (\sid.9)  we can
give a
necessary and sufficient numerical criteria for this to happen.

\cor 13 A hyperplane $H\subset(\C^{n+1},0)$ is a limiting tangent hyperplane of
the hypersurface $X=V(f)$ at 0 if, and only if, the numbers $e_i(J(f)_H)$ and
$e_i(J(f))$ differ for some $i=1,\dots,n.$
\endgroup\demobox


\sectionhead {\af} Strict Dependence and Thom's Condition $a_f$


We review the notion of strict dependence and characterize this condition in
terms of vanishing of pullbacks along components of exceptional divisors of
some normalized blowups. This allows us to prove a
principle of specialization of strict dependence in (\af.3).

Strict dependence is used to describe equisingularity conditions like
Whitney's condition condition $a$ and Thom's
condition $a_f$ in terms of integral closure; see \cite{\Ma} for a discussion
of Thom's condition $a_f.$ We apply our results to study
this condition for a function
$f$ on $(\C^{n+1},0)\times (\C^p,0).$ This may be viewed as a family of
functions on $(\C^{n+1},0).$  Proposition (\af.5) recovers a result
of  Massey. It gives sufficient numerical conditions for $a_f$ to hold.
The same approach gives a similar result for a map $f$ of reduced analytic
spaces. It involves the relative Jacobian ideal and the relative Nash blowup
the definition of which we review.

\stp1 In the Setup (\sid.1), we say that an element $h$ of $\OO_{X,0}$ is
strictly dependent on an ideal
$I$ in $\OO _{X,0}$ if for all holomorphic curves $\phi:(\C,0)\to(X,x)$,
the pullback $\phi^*h$ is contained in the ideal  $m_1\phi^*(I)$ where $m_1$ is
the maximal ideal in
$\OO_{\C,0}$.  All such elements form the ideal $\bar I^\dag.$

We will consider
the normalized blowups $B_{2,N}$ and $B_N$ of $X$ along $I$ and $m_YI$ with
structure maps
$b_{2,N},b_N$ and
exceptional divisors $D_{2,N}$ and $D_N.$ We denote the pullback to $B_N$ of
the divisor $D_2$ by $\tD{2,N}.$

\prop2 Let $h$ be a holomorphic function on $X.$ Then, the following
conditions are equivalent.

\itemitem{(i)} For every point $y$ in $0\times Y,$ the germ of
$h$ at $y$ is strictly dependent on $I\OO_{X,y}.$
\itemitem{(ii)} For every point $z$ in $B_{2,N}$ over a point in $Y$ and any
element $g$ of $I$ the pullback of which to $B_{2,N}$ generates a local
equation of $D_{2,N}$ at
$z$ the quotient $h\circ b_{2,N}/g\circ b_{2,N}$ is a holomorphic
function near $z$ and vanishes on
$|D_{2,N}|\cap b_{2,N}\inv(Y)$ near $z.$
\itemitem{(iii)} For each component $C$ of the underlying set of $\tD{2,N}$ the
order of vanishing of $h\circ b_N$ along $C$ is greater than the order of
vanishing of the pullback $b_N^*I$ along $C$ if $b_N(C)$ is contained in
$0\times Y,$ and, if $C$ is not mapped to $Y,$ the order of vanishing of
$h\circ b_N$ along $C$ is not smaller than the order of vanishing of the
pulback of $I$ along $C.$

\pf Assume that  condition (i) obtains. Then, in particular, for a point
$y$ in$Y$ the germ of
$h$ at $y$ is in the integral closure of
$I\OO_{X,y}.$ Hence, for each point $z$
in $|D_{2,N}|$ over $Y$ and an element $g$ as in (ii) there exists a
neighborhood $U$ of $z$ and a holomorphic function $k$ on $U$ such that
$h\circ b_{2,N}=k(g\circ b_{2,N})$ on $U$. We want to show that $k$
vanishes on the fibre of $|D_{2,N}|\cap U$ over $Y.$ So, pick a point
$z'$ in this fibre close to $z.$ Then, for a curve $\phi_N:(\C,\C-0,0)\to
(B_{2,N},B_{2,N}-|D_{2,N}|,z')$ we have
$$h\circ b_{2,N}\circ \phi_N=(k\circ\phi_N)
(g\circ b_{2,N}\circ\phi_N)\in m_1 {\phi}^*I.$$ Here $\phi$ is the path
$b_{2,N}\circ\phi_N.$ This implies that $k$ vanishes at $z'$.

On the other hand, if (ii) holds, then, clearly, the function $h$ vanishes on
$V(I).$ Hence it is enough to consider  paths
$\phi:(\C,\C-0,0)\to(X,X-V(I),y).$ Now, $\phi$ can be lifted to a path
$\phi_N$ on the normalized blowup
$B_{2,N}$. Denote the image of 0 under the lifted path by $z.$ Then, there is
neighborhood $U$ of $z$ in $B_{2,N}$ and a holomorphic function $k$ on $U,$
vanishing on $|D_{2,N}|\cap b_{2,N}\inv(Y)\cap U,$ so that
$h\circ b_{2,N}=k(g\circ b_{2,N})$ where $g$ is as in the statement of (ii). In
particular, the pullback $k\circ\phi_N$ vanishes at z. Hence,
$h\circ b_{2,N}\circ \phi_N=h\circ\phi \in m_1\phi ^*(I)$.

Finally, condition (ii) and (iii) are equivalent by an argument similar to the
ones already used. For a point $z$ in $B_{2,N}$ over a point of $0\times Y,$
we have $h\circ b_{2,N}=k(g\circ b_{2,N})$ with $k$ and $g$ as above.Then
$k$ vanishes on the fibre over $0\times Y$ of the underlying set of
$D_{2,N}$ near $z$ if, and only if, the pullback of
$k$ to $B_N$ vanishes on the part of $\tD{2,N}$ formed by components that are
mapped to $0\times Y.$ This finishes the proof.\enddemo

\thm3 In the setup of (\sid.1), let
$h\in\OO_{X,0}$ be a function so that $h\OO_{X(t),0}$ is strictly dependent
on $I(t)$ for all $t$ in a Zariski--open subset of $Y.$ If the numbers
$e_1(t),\dots,e_n(t)$ are independent of $t,$ then $h$ is strictly
dependent on $I.$

\pf The constancy of the $e_i(t)$ implies that for each component  of
$D_2$ all components of its intersection with $b_2\inv(Y)$ map onto $Y$; hence
the same statement is true for the normalized blowup.

Now, we use an argument similar to the one used in the above proof. Let $z$
be a point in $B_{2,N}$ over 0 in $X.$ From the principle of specialization,
we know that $h$ is contained in the integral closure of $I.$ Hence, near
$z,$ we can write $h\circ b_{2,N}=k(g\circ b_{2,N})$ with $g$ and $k$ as in
the above proof. Let $C$ be component of the underlying set of the
exceptional divisor that contains $z.$ Then, its intersection with the
inverse image of $Y$ maps onto $Y.$ In particular, the set of points $z'$
in this intersection that map to the Zariski open subset of $Y$ where
the strict dependence condition obtains are dense. By (\af.2), the
holomorphic function $k$ vanishes at such a point $z'.$  Therefore,
the function $k$ vanishes on the intersection. Hence, again by
(\af.2), the claim follows.\enddemo

\art 4 The case of a smooth ambient space

Let $M=(\C^{n+1},0)\times(\C^p,0)$ and consider a function germ
$f:M\to(\C,0).$ Let $Y=0\times(\C^p,0)$ and $\Sigma=V(J(f))$ be the critial
locus of
$f.$ We assume that $\Sigma$ is nowhere dense and doesn't contain
$(\C^{n+1},0)\times 0.$ Let $z_0,\dots,z_n$ be coordinates on $(\C^{n+1},0).$
For a series of points $p_i$ in $X-\Sigma$ converging to 0, we may view the
series of tangent planes to the level hypersurfaces of $f$ through the $p_i$
as a series of points in $\PP^n.$ We may assume that this series converges
after passing to a subseries of $\{p_i\}.$ We call this limit a {\it limiting
tangent hyperplane} to the fibres of $f$ at 0. If all limiting tangent
hyperplanes contain $Y,$ we say that the pair $(M-\Sigma,Y)$ satisfies $a_f$
at the origin. This condition translates into the integral dependence condition
$${\pd f\over\pd v}\in\overline{\left({\pd f\over\pd z_0},\dots,{\pd f\over\pd
z_n}\right)\OO_M}^\dag\text{ for all }v\in T_0Y.$$ (See \cite{G2} for this and
the definition of the strict integral closure, denoted by the superscript
$\dag.$) It is known that this condition obtains for generic points of $Y.$
Hence we can apply the principle of specialization (\sid.7), which also works
for the strict closure. For a point $t$ in $Y$ we denote
$$\lambda_i(t):=e_i\left(\left({\pd f\over\pd z_0},\dots,{\pd f\over\pd
z_n}\right),(\C^{n+1},0)\times  t\right).$$ Then we get the following
result of  Massey \cite{\Ma}.

\prop 5 In the above setup the  pair $(M-\Sigma,Y)$ satisifes $a_f$ at
the origin if the numbers $\lambda_i(t)$ are independent of $t$ in $Y.$
\endgroup\demobox

\art 6 Relative Jacobian ideal and relative Nash modification

We recall some definition and facts from \cite{T4,Ch.II} and \cite{T2,
Chapter 2}. Suppose $f:(X,0)\to(S,0)$ is a morphism of analytic germs so that
the sheaf of relative differentials $\Omega^1_f$ is locally free of rank
$d=\dim X-\dim S$ outside a nowhere dense analytically closed subset $F$ of
$X.$ We will consider the {\it relative Jacobian ideal} $J(f).$ It is defined
as the $d$th Fitting ideal of the $\OO_X$--module $\Omega^1_f.$

Alternatively, it can be defined via an embedding $X\hookrightarrow
S\times(\C^N,0)$ so that $f$ is the restriction to $X$ of the projection onto
the second component. If $X$ is defined by $f_1,\dots,f_m$ in
$S\times(\C^N,0)$ and $z_1,\dots,z_N$ are coordinates of $\C^N,$ then $J(f)$
equals the ideal of $\OO_{X,0}$ generated by the minors of rank $N-d$ of the
{\it relative Jacobian matrix } $$\left({\pd f_i\over \pd
z_j}\right),i=1,\dots,m;j=1,\dots,N.$$

On the other hand, we may consider the Grassmannian space
$\text{Grass}_d\Omega^1_f\to X$ of locally free quotients of $\Omega^1_f$ of
rank $d.$ There is a well--defined section over $X-F.$ The closure of its
image is called the {\it relative Nash modification}. Geometrically, this
section maps a point $p$ to the tangent space at $p$ of the level surface of
$f$ passing through $p.$

If $X$ is a relative complete intersection over $S,$ the relative Nash
modification equals the blowup of the relative Jacobian ideal in $X.$ In the
general case we have to embed $X$ into a relative complete intersection $X'$
over $S$ and and blow up the ideal $J(f)'$ in $\OO_X$ induced by the
relative Jacobian ideal of $X'$ over $S.$

\art 7 The general case

Let $X\sub(\C^N,0)\times(\C^p,0)$ be a reduced analytic space germ of
pure dimension $n+p,$ so that the fibres of the map $X\to Y=(\C^p,0)$
induced by the projection onto the second factor  form a
family of reduced analytic germs of pure dimension $n.$ The fibre over
$t\in Y$ will be denoted by $X(t).$ We assume that $0\times(\C^p,0)$ is
contained in $X$ and will identify it with $Y.$ Let $(S,0)$ be a reduced
analytic space germ and $f:X\to(S,0)$ be an analytic map germ that maps $Y$ to
0 and is a submersion off a nowhere dense analytic subset $F$ of $X.$ The $a_f$
condition is defined as above. It is known that it is satisfied at generic
points of $Y,$ if there is no blowup in codimension 0 (see \cite{\HMS}). Then,
$(X-F,Y)$ satisfies
$a_f$ at 0 if each component of the exceptional divisor of the relative Nash
blowup maps onto
$(\C^p,0).$ The principle of specialization (\sid.7) gives sufficient numerical
conditions for this to happen.

\thm 8 In the above situation (\af.4) the pair $(X-F,Y)$ satisfies $a_f$ at 0
if the Segre numbers
$e_1(J(f)',X(t)),\dots,e_n(J(f)',X(t))$ are independent of $t$ in $Y.$
\endgroup\demobox


\sectionhead{\wf} Condition $W_f$ and the Whitney conditions


We will now use the Theorem of Specialization of Integral Dependence to
study families of functions on $(\C^{n+1},0).$ Theorem (\wf.2) gives a
sufficient criterion for the ambient space to satisfy the $W_f$ condition
along a the parameter space.   It requires the constancy of numbers that
only depend on the
family members. In particular, the constancy of these numbers  implies that
the smooth part of the corresponding family of hypersurfaces is
Whitney regular along the parameter axis.

Our condition is also necessary for this family of hypersurfaces to be
Whitney equisingular along the parameter space. One part of our numbers can
be computed from the Milnor fibres of the restriction of the family members to
generic linear subspaces.

We also apply our theory to the study of families of 2-dimensional
hypersurfaces (\wf.6), families of hyperplane slices
of a hypersurface (\wf.7), and the deformation of a hypersurface to its
tangent cone (\wf.8).

\stp 1 Consider the family of smooth spaces
$M=(\C^{n+1},0)\times(\C^p,0)\buildrel F\over\to(\C^p,0)=Y,$ and identify $Y$
with the zero section $0\times(\C^p,0).$ Let $f:M\to(\C,0)$ be a map germ that
maps
$Y$ to 0. Let $I$ be the Jacobian ideal of $f,$ generated by its partial
derivatives. Denote the restriction of $f$ to a fibre $M(t)$ by $f_t.$ We
assume that $f_t$ is reduced for any
$t.$ Also, assume that the critical locus
$\Sigma$ of $f$ contains
$Y.$ Fix coordinates
$z_0,\dots,z_n,t_1,\dots,t_p$ which respect the product structure of $M.$ We
will use the relative Jacobian ideal $J_z(f),$ generated by the partial
derivatives of $f$ with respect to $z$--coordinates.

Using  a suitable
distance function `dist' for linear subspaces of $M,$
condition $W_f$ can be expressed as an inequality; see e.g.
\cite{GK, (4.1)}.  The pair $(M-\Sigma, Y)$
satisfies
$W_f$ at 0 if there exists a neighborhood $U$ in $M$ of 0 and a postive
constant $C$ so that for all points $x$ in $U-\Sigma$ we have
$$\dist(T_0Y,T_xV(f-f(x)))\leq C\dist(x,Y),$$
Using (1.1)(ii), one can show that this condition is equivalent to the
following relations
$${\pd f\over\pd t_i}\in\overline{m_zJ_z(f)\OO_{\C^{n+1+p},0}},\quad\text{ for
} i=1,\dots,p,\tgs1.1$$ where $m_z=(z_0,\dots,z_n)$ is the ideal defining $Y$
in
$M.$ The latter condition was named (c)--equisingularity by Teissier; see
\cite{T2. 2.17, p.601}.

Consider now the family $X$ of hypersurfaces defined by
$f\in\OO_{\C^{n+1+p},0}.$  We use the above distance function to express the
Whitney conditions as an inequality also. The pair
$(X-\Sigma, Y)$ satisfies the Whitney conditions at 0 if there exists a
neighborhood $U$ in
$X$ of 0 and a postive constant $C$ so that for all points $x$ in $U-\Sigma$
we have
$$\dist(T_0Y,T_xX)\leq C\dist(x,Y).$$ We also say that $X-\Sigma$ is
Whitney regular along $Y$ at 0.

Again making use of 1.1 (ii), this condition is equivalent to the following
integral dependence relation in the local ring $\OO_{X,0}$ of $X$ at 0 (see
\cite{G1}):
$$\left.{\pd f\over\pd t_i}
\right|_X\in\overline{m_zJ_z(f)\OO_{X,0}},\quad\text{ for }
i=1,\dots,p.\tgs1.2$$

Note that the $W_f$ condition implies the
Whitney conditions as an integral dependence relation descends to quotient
rings.

\medskip
More generally, Whitney regularity is defined for any pair of submanifolds of
$M$ in an analogous way. A Whitney stratification of a singular space is a
stratification such that for any strata $S$ and $S'$ with $S'$ contained in the
closure of $S$ the bigger strata $S$ is Whitney regular along $S'.$ We say that
the hypersurface
$X$ is {\it Whitney equisingular} along $Y$ if there exists a Whitney
stratification of $X$ with $Y$ as a stratum.

\medskip
If $X$ is Whitney equisingular along $Y,$ then $X$ is topologically trivial
along $Y$ by the Thom--Mather isotopy theorem; see e.g. \cite{T4, Ch. VI,
4.3.1, p.482}. This trivialization can be obtained by lifting vector fields on
$Y$ to a corrugated (Fr. {\it rugueux}) vector field on $X.$ These vector
fields are integrable and are tangent to the strata of $X;$ see \cite{V}.

If, in addition, $(M-\Sigma,Y)$ satisfies the condition $W_f,$ any vector
field on $Y$ can be lifted to a corrugated vector field on the ambiant space
$M$ that is tangent to the level hypersurfaces of $f;$ see \cite{T2, 2.17,
Cor. 1, p.602}. The integration of these vector fields yields a
right--trivialization of $F.$ In particular, the maps $f_t$ are
equivalent up to homeomorphism of their source $(\C^{n+1},0).$

\medskip
For $k=2,\dots,n+1,$ we use the following notations:
$$\align
\lambda_k(f_t)&=e_k(J_z(f),M(t))=e_k(J(f_t),M(t)),\\
m_{k-1}(f_t)&=m_{k-1}(J_z(f),M(t))=m_{k-1}(J(f_t),M(t)).\endalign$$ The number
$\lambda_k(f_t)$ is called the
$k$th  {\it L\^e number} of $f_t.$ Note that
$e_1(J_z(f),M(t))=e_1(J(f_t),M(t))$ is zero as, by assumption, the critical
locus of $f_t$ has codimension at least two in $M(t).$ The number $m_k(f_t)$
is called the $k$th {\it relative polar multiplicity of $f_t.$}

\medskip
Fix $t$ in $Y.$ For $L^k\sub(\C^{n+1},0)$ a generic $k$--dimensional linear
subspace, we denote the Euler characteristic of the Milnor fibre of
${f_t}|_{L^k}$ at 0 by
$\chi^{(k)}(t).$ Let $s$ be the codimension of the singular locus of $f_t.$
Then, by a result of  Massey \cite{Ma, 10.6, p.96} and (\form.4.2), we
have
$$\eqalign{
\chi^{(k)}(t)&=(-1)^k\sum_{i=2}^k(-1)^i(e_i(J_z(f),L^k\times t)\cr
&=m_k(f_t)+\lambda_k(f_t)-\lambda_{k-1}(f_t)+\lambda_{k-2}(f_t)-\dots\pm
\lambda_s(f_t).}$$ Hence, the number $\chi^{(k)}(t)$ doesn't
depend on the choice of the
$L^k.$ We define
$$\chi^*(t):=(\chi^{(n+1)}(t),\dots,\chi^{(2)}(t)).$$

\thm2 Suppose that the map
$t\mapsto (m_1(f_t),\dots,m_n(f_t),\chi^*(t))$ is constant on $Y.$
Then the pair $(M-\Sigma,Y)$ satisfies $W_f$ at 0. In particular, the smooth
part of $X$ is Whitney--regular along $Y$ at 0. Also, the smooth parts
 of
$\Sigma$ and the components of the singular locus of $\Sigma$ of codimension
one in $\Sigma$ satisfy the Whitney--conditions along
$Y$ at 0.

\pf We now use the strength of the machinery developed in this
paper. The constancy of the map in the proposition is equivalent to
the constancy of the map
$$t\mapsto
(m_1(f_t),\dots,m_n(f_t),\lambda_2(f_t),\dots,\lambda_{n+1}(f_t)).$$ Hence,
by the expansion formula (\form.5), the numbers $e_k((m_zJ_z(f))(t),M(t))$ are
independent of $t$ in $Y.$ Also, we know that the integral dependence
relation (\wf.1.1) is satisfied at generic points of $Y.$ Therefore, we can
apply the Principle of specialization of integral dependence (\sid.7) to
prove the first assertion.

The second assertion follows from the following observation:

Consider a
stratum $W$ in $\Sigma$  whose closure  is the image of a component $C$ of the
exceptional divisor $D$ of the blowup of $M$ along $J(f).$ Denote its
structure map by $b.$ Then the conormal of
$W$ in $M$ equals $C.$ In fact, it is well--known
that the pair
$(M-\Sigma,W)$ satisfies
Thom's $a_f$ condition at points of a Zariski--open dense subset $U$ of $W.$
Hence, we have $$C\cap b\inv(U)\sub C(W)\cap b\inv(U).$$ Replace $U$ be a
perhaps smaller Zariski--open dense subset of $W$ so that for $w\in U,$ the
fibre over $w$ of $C(W)$ is isomorphic to $\PP^{k-1}$ where $k$ is the
codimension of $W.$ A dimension count shows that, after shrinking $U$ once
more, we may assume that the fibre of $C$ over a point $w$ in $U$ has dimension
$k-1.$ It follows that the fibres of $C$ and
$C(W)$ over $w$ are equal. As $C(W)$ is the closure of $C(W)\cap b\inv(U)$ in
$H\times\PP(H),$ the claim follows. (The claim also follows from the
Principle of Lagrangian Specialization \cite{LT1, 1.2.6}.)

We are going to show now that the smooth part of $W$ is Whitney regular along
$Y$ if the hypothesis of the theorem obtains. By Teissier's
characterization of Whitney--conditions \cite{T4,
Ch.5}, we have to show that the exceptional divisor of the blowup of $C$
along the preimage of $Y$ in $C$ is equidimensional over $Y.$ As (\wf.1.1)
holds, the ideal $J_z(f)$ is a reduction of $J(f);$ thus the
component $C$ is finite over over a component $C'$ of the exceptional divisor
of the blowup of $M$ along $J_z(f).$  Hence, it is enough to show the analogous
statement for $C'.$ But this is part of the result of (\sid.6).

Clearly, a component of $\Sigma$ is the image of a component $C$ of $D.$ Also,
every component of the singular locus of $\Sigma$ of codimension one in
$\Sigma$ is the image of  a component $C$ of $D$ (see \cite{\Ma, Prop. 1.32,
p.30}). This finishes the proof.\enddemo

\thm3 Suppose that $X$
admits a Whitney stratification with $Y$ as a stratum. Then, the map
$t\mapsto (m_1(f_t),\dots,m_n(f_t),\chi^*(t))$ is constant on $Y.$

\pf
The proof uses topological methods. By the Mather--Thom isotopy theorem, the
topological type of $X(t)$ is independent of $t$ in $Y$ (see e.g. \cite{T4,
Ch. VI, 4.3.1}). Hence, by an observation of L\^e D\~ung Tr\'ang \cite{L, p.
261}, the homotopy type of the Milnor fibre of $f_t$ is independent of $t.$ In
particular, the Euler characteristic is constant. That is
$\chi^{(n+1)}(0)=\chi^{(n+1)}(t).$

Fix a point $t$ in $Y.$ For $k=2,\dots,n,$ we can choose a $k$--dimensional
linear subspace
$L^k$ of
$(\C^{n+1},0),$ generic w.r.t $f_0$ and $f_t,$ so that the restriction of $f$
to $L^k\times Y$ satisfies the same assumptions as $f.$ Hence, we get by the
same argument as above $\chi^{(k)}(t)=\chi^{(k)}(0).$ It follows that
$\chi^*(t)$ is independent of $t$ in $Y.$

It remains to show the constancy of the map
$t\mapsto(m_1(f_t),\dots,m_n(f_t)).$ We prove this by induction over the
dimension of the fibres of $X$ over $Y.$ The case of a one--dimensional
family follows from the classical theory of equisingularity for plane curves;
Whitney equisingularity implies that both, $\lambda_2(f_t)=\mu^{(2)}(X(t))$
and $m_1(f_t)$ are independent of $t$ in $Y.$

Assume now that the theorem has been proven for fibre dimension $n-1.$ Fix a
point $t$ in $Y.$ Then, we may choose a hyperplane $H$ in $M(0)$ so that
$X\cap(H\times Y)$ is again Whitney--equisingular along $Y,$ and $H$ is
generic w.r.t $X(t)$ and $X(0).$ In particular, we may assume that for
$k=1,\dots,n-1$
$$m_k(f_t)=m_k(J_z(f),H\times t)=m_k(J_z(f),H\times 0)=m_k(f_0)$$ This
implies, together with the constancy of $\chi^{(*)}(t)$, that
$$t\mapsto(m_1(f_t),\dots,m_{n-1}(f_t),
\lambda_2(f_t),\dots,\lambda_{n-1}(f_t) \tgs3.1$$
is constant on $Y.$ Therefore, using the expression of $\chi^{(n)}(t)$ in terms
of
$m_k(f_t)$ and $\lambda_k(f_t),$ we see that the map $t\mapsto
m_n(f_t)+\lambda_n(f_t)$ is constant.

Furthermore, by (\wf.3.1) and (\sid.3), the polar
variety
$P_n(J_z(f),H\times Y)$ specializes; therfore, it is empty. It follows that
every component of
$P_n(J_z(f),M)$ contains
$Y.$ Hence, it remains to show that the fibre over 0 of a component of
this polar variety is not contained in the fibre of $\Sigma$ over 0. Suppose
there were such a component. Then, there exists a component $C$ of the
exceptional divisor of the blowup of $J_z(f)$ in $M$ that doesn't map onto
the parameter space $Y$ and so that $C(0)$ maps to a subset of
$M(0)$ of non--zero dimension. Consider a point
$p$ in this set close to 0, the stratum $S$ of the initial
Whitney--stratification through
$p,$ and the familiy $X\to S$ near $p$ given by some retraction. The fibres of
this family are of dimension at most $n-1.$ Hence, by the induction hypothesis
and (\sid.6), every component of the exceptional divisor of the
blowup of
$J_z(f)$ in the germ of $M$ at $p$ maps onto $S.$ This contradicts the
existence of $C.$ It follows that $P_n(J_z(f),M)$ specializes. Therefore, the
map $t\mapsto m_n(f_t)$ is upper semi--continuous, but so is $t\mapsto
\lambda_n(f_t)$ by the lexicographically upper
semi--continuity of the L\^e numbers (\sid.5). We know that the sum of the
two maps is constant, and therefore each one is. This finishes the proof.
\enddemo

\cor4 Assume that $X$ admits a Whitney stratification with $Y$ as
a stratum, then $M-\Sigma$ satisfies the $W_f$ condition along
each stratum.

\pf For each stratum $S,$ apply Theorem (\wf.3) to the family $X\to S$ given by
some retraction which is compatible with $X\to Y.$  Then, apply (\wf.2).
\enddemo

\rmk5 (1) This corollary recovers a result of  Parusinski
\cite{Pa} which was also proven by Briancon, Maisonobe and Merle \cite{BMM} in a
more general context.

In his proof Parusinski shows that the assumptions of the corollary imply that
each component $C$ of the exceptional divisor of the blowup of $J(f)$ in $M$ is
the conormal of the closure of some stratum $S$ of the Whitney
stratification of $X.$ By assumption, the smooth part of $S$ is Whitney
regular along $Y.$ Hence, the exceptional divisor of the blowup of
$C$ along the preimage of $Y$ is equidimensional over $Y.$ Thus, by (\sid.6),
the Segre numbers of $J(f)$ on the fibres $M(t)$ are
independent of $t.$ Furthermore, the $W_f$ condition implies that $J_z(f)$ is a
reduction of $J(f).$ So, these Segre numbers equal $\lambda_k(f_t).$ Also, it
is not hard to see that the
$W_f$ condition implies the constancy of
$m_k(f_t).$ Therefore, we can use Parusinki's proof to give a purely
algebro--geometric proof of the above Theorem (\wf.3). (Note that all results
of  L\^e and Teissier's work \cite{LT1} can also be proven by purely
algebro--geometric means. This was worked out by Roberto Callejas--Bedregal
in his thesis \cite{Ca}).

\smallskip
(2) Consider the diagram of blowups in (\sid.1) for $X=M,I=J_z(f)$ and $m$
replaced by $m_z.$ Results of Henry, Merle and Sabbah \cite{HMS, 6.1, 3.3.1}
show that $(M-\Sigma,Y)$ satisfies $W_f$ at 0 if, and only if, the
exceptional divisor $\tD1$ is equidimensional over $Y.$

On the other hand, the constancy of $t\mapsto
(m_1(f_t),\dots,m_n(f_t),\chi^*(t))$ implies that the L\^e numbers
$\lambda_k(f_t)$ are constant along $Y.$ Hence, a necessary condition for
the constancy of this map is that every component of
$\tD2^{X-Y}\cdot\tD1$ and $D_2^Y$ is equidimensional over $Y.$ As we have
seen in remark (1), this condition controls some lower--dimensional strata of
the Whitney--stratification of $X.$

\smallskip
(3) Denote the polar multiplicities, resp. Segre numbers of $J(f_t)$ on $X(t)$
by $m_k(X(t)),$ resp. $\lambda_k(X(t)).$
Then, the principle of specialization and the expansion formula show that the
constancy of
$$t\mapsto(m_1(X(t)),\dots,m_{n-1}(X(t)),
\lambda_1(X(t)),\dots,\lambda_n(X(t)) )$$ implies that the pair
$(X-\Sigma,Y)$ satisfies the Whitney--conditions at 0. However, to get a
converse statement, one has to control the components of the exceptional
divisor $D$ of the blowup $B$ of $X$ along the ideal induced by $J_z(f).$  It
seems that the existence of a Whitney stratification of $X$ with $Y$ as a
stratum is not sufficient for the constancy of this map; for example, suppose
$X$ is a 1--parameter family  of non--isolated surface singularities in
$(\C^3,0).$ Denote the structure map of $B$ by $b.$  Suppose the exceptional
divisor $D$ has two components
$C$ and
$C'$, where the image of
$C'$ is the parameter stratum $Y$, and the image of $C$ is a component $W$  of
the singular set of the total space. Suppose $C'$ is equidimensional over
$Y$,  and contains $C\cap b\inv(Y)$.  It is conceivable that
$C\cap b\inv(Y)$ is not equidimensional over $Y$, even if the smooth part of
$W$ is Whitney regular over $Y$, and the smooth part of $X$ is Whitney regular
over $Y$.
 In fact, Whitney regularity of the smooth part of $X$ along $Y$ implies that
the exceptional divisor $D_Y$ of the blowup of $B$ along $b\inv(Y)$  is
equidimensional over $Y$ (\cite{T4, Ch. 5, 2.1, p.470}).  However, since
$C\cap b ^{-1}(Y)$ is hiding inside the hypersurface $D$,  there need not
be a component of $D_Y$ over this intersection; so Whitney regularity of the
smooth part along the parameter stratum may not control the behavior of
$C$ over $Y$.   If the smooth part of $W$ is Whitney regular  over $Y$, then,
as we have seen in the proof of (\wf.3), we can control the behavior of
$C(W)$, the conormal of $W$. However,
$C$ is a proper subset of $C(W)$; generically it consists of smooth points of
$W$ and tangent
 hyperplanes that are limits of tangent hyperplanes from the smooth part of
$X$.  So it may not be controlled by this condition either.

This hypothetical example underlines the possible geometric difference between
the Whitney conditions and the constancy of the above numbers. In both
conditions  the components of $D$ that map onto $Y$ are equidimensional over
$Y$. However, the constancy of the above Segre numbers requires that the
preimage of
$Y$ in each component of
$D$ be equidimensional over $Y$ while the Whitney conditions imply
the equidimensionality over $Y$ of the exceptional divisor $D_Y$ of the blowup
of
$B$ along
$Y$.  It is an open question whether or not  the above
hypothetical example exists or not.

\medskip
The conclusions of Theorem (\wf.2) don't give us any information on
whether the smooth part of $X$ is Whitney regular along the smooth part of
$\Sigma,$ or the singular locus of $\Sigma.$ Therefore, in general, we can't
use (\wf.2) to show that $X$ is Whitney equisingular along $Y.$ However, in
specific situations one may use auxiliary information to build a Whitney
stratification with $Y$ as a stratum.

\cor6 Assume that $X$ is a one--parameter family of reduced non--isolated
surface singularities. Then, $X$ is Whitney--equisingular along the parameter
axis $Y$ if and only if the map
$t\mapsto(m_1(f_t),m_2(f_t),\lambda_2(f_t),\lambda_3(f_t))$ is constant on
$Y.$

\pf The `only if' statement follows from Theorem (\wf.3). So, assume
now that the map is constant. Thus, by Theorem (\wf.2), the smooth parts of
$X$ and its singular locus $\Sigma$ are Whitney regular along $Y.$ It
follows that $\Sigma$ is either smooth, or its singular locus contains
$Y.$ Assume for the moment that
$\Sigma$ is singular. Then, its singular locus is 1--dimensional, hence of
codimension one in $\Sigma.$ Therefore, by (\wf.2), the singular locus of
$\Sigma$ equals $Y.$ The remaining step to showing that
$$(X-\Sigma,\Sigma-Y,Y)$$ is a Whitney stratification is to establish the
Whitney regularity of $X-\Sigma$ along $\Sigma-Y.$ To see this, consider a
general linear form  $l$ on $\C^3$ the kernel of which intersects the
fibres $X(t)$ and
$\Sigma(t)$ transversally for all $t\in (Y,0).$ Also, assume that for any
any point $x\in X-Y$ the kernel of $l$ is not a limiting tangent plane
of $X(t)$ at $x.$ A form satisfying the transversality conditions exists
because the smooth parts of $X$ and $\Sigma$ are Whitney regular along
$Y.$ Also, by (\sid.6)(ii), the fibre of any component of $D^{X-Y}$ over 0 is
of dimension at most 2. Here, $D$ denotes the exceptional divisor
of the blowup of  $M$ along $J_z(f).$ Hence, we can choose the form $l$ so
that the ideal $J(f)_{z,l}$ generated by the partial derivatives of $f$ by
vectors in the kernel of $l$ is a reduction of $J_z(f)$ outside $Y.$ By
(\sid.11), this implies the desired behavior with respect to limiting tangent
planes.

Next, construct the family of reduced plane curves
$$\pi:X-Y\to Y\times\C,(x,t)\mapsto(t,l(x)).$$ The Whitney--equisingularity
of this family will imply the desired result. But, for a family of plane curves
this is equivalent to the constancy of the Milnor number of its fibres
$X(t,s).$ Let
$f$ be a defining equation of $X\subset(\C^3\times\C,0).$  Then, using
(\form.4.2), we have
$$\mu(X(t,s))=\sum_{x\in\Sigma(t,s)}\lambda_2(f_{t,s},x)=
\sum_{x\in\Sigma(t,s)}\lambda_2(f_{t},x)+\sum_{x\in\Sigma(t,s)}m_2(f_{t},x).$$
Here we have used that $J(f)_{z,l}$ is a reduction of $J_z(f)$ outside $Y.$
So, $\lambda_2(f_{t,s},x)=e_2(J_z(f),(M(t,s),x)).$

Now, as the smooth part of $\Sigma$ is Whitney equisingular along $Y,$ the
first  sum equals $\lambda_2(f_t)$ by (\segre.3). (The equisingularity
ensures that the result of (\segre.3) remains valid for arbitrary
$s.$)  Furthermore, we claim that the second sum vanishes.
Otherwise, there would be a fixed component of
$\Lambda_3(j_z(f),\C^3\times\C)$ outside of
$Y.$ This is only possible if a component of the exceptional divisor
of the blowup of $\C^3\times\C$ along $J_z(f)$ maps to this component of
the L\^e cycle
But, by (\sid.6), this cannot occur. This finishes the proof.\enddemo

\expl 7 We can apply this corollary to families of hyperplane sections of
hypersurfaces of dimension 3. Let $\Cal X=V(f)\subset(\C^4,0)$ be a reduced
hypersurface with singular locus $\Cal S$ of dimension two. Let $H_0$ be a
hyperplane in $(\C^4,0)$ so that the intersection $\Cal X\cap H_0$ is reduced.
Suppose $H_0$ is a regular point of a curve $V$ in the projective
space $\PP^3$ of all hyperplanes in $\C^4$ through 0. Then  the family of
reduced surfaces
$$X=\{\Cal X\cap H\}_{H\in(V,H_0)}\subset (\C^3\times V,(0,H_0))$$ is Whitney
equisingular along $Y=0\times (V,H_0)$ if, and only if, the map
$$H\mapsto(m_1(f|_H),m_2(f|_H),\lambda_2(f|_H),\lambda_3(f|_H))$$ is constant
on $(V,H_0).$

\medskip
If we assume an extra condition, then we get the following stronger statement.
\smallskip
{\it In the above setup, assume,in addition, that $H_0$ intersects $\Cal S$
and the singular locus of
$\Cal S$ transversally. Then, the
family of reduced hypersurfaces
$$X=\{\Cal X\cap H\}_{H\in(\PP^3,H_0)}\subset (\C^3\times \PP^3,(0,H_0))$$ is
Whitney equisingular along $Y=0\times (\PP^3,H_0)$ if, and only if, the map
$$H\mapsto(\lambda_2(f|_H),\lambda_3(f|_H))$$ is constant
on $(\PP^3,H_0).$}

\smallskip
The additional transversality conditions is needed to ensure that the
dimension of any component of the singular locus of the singular locus
$\Sigma$ of $X$ is of dimension 3. Hence, the same arguments as above yield
the result, provided that the constancy of the L\^e numbers implies that the
map $H\mapsto(m_1(f|_H),m_2(f|_H))$ is constant.

To see this, we use the observation \cite{G3, 2.9} of the first author. Let
$g$ be a defining function of the hypersurface $X$ in
$M=(\C^3\times\PP^3,(0,H_0)).$ Then, $(M-\Sigma,Y)$ satisfies $w_g$ at
$(0,H_0)$ if and only if the map $H\mapsto(\lambda_2(g|_H),\lambda_3(g|_H))$
is constant on $Y.$ Now, fix a point $(0,H_1)$ on $Y$ and let $L$ be the line
joining $H_0$ to $H_1.$ Denote the restriction of $g$ to $\C^3\times L$
by $g'.$ Then,
$(\C^3\times L-\Sigma,0\times L)$ satisfies
$w_{g'}$ at $(0,H_0)$ as the integral dependence relation $(\wf.1.1)$
restricts to the desired relation on the subspace in question. In particular,
if
$(z_1,\dots,z_3)$ are coordinates on $\C^3,$ the ideal $J_z(g')$ is a
reduction of $J(g').$ Hence, by \cite{\HMS, Thm. 6.1}, the polar multiplicities
of $J_z(g')$ on $(\C^3\times L,(0,H))$ are independent of $H$ in $L.$
Furthermore, by a polar transversality theorem \cite{HMS, 4.4.6} and
(\sid.3), these polar multiplicities equal the relative polar multiplicities of
$g|_H.$ Hence, the map $H\mapsto m_1(g|_H),m_2(g|_H))$ is constant on
$Y.$ This finishes the proof of the claim.

\art8 The deformation to the tangent cone

Let $(X,0)=V(f)\subset(\C^{n+1},0)$ be a
hypersurface germ. Then, recall that the deformation of $X$ to its tangent
cone $\Cal X\to \C$ is defined by $g=t^{-d}f(tz_0,\dots,tz_n),$ where
$d$ is the degree of the initial term $f_0$ of $f.$ The fibre of $\Cal
X$ over 0 equals the tangent cone $C_0X=V(f_0)$ of $X.$ For non--zero
$t,$ the fiber $\Cal X(t)$ is isomorphic to $X.$ The parameter space is
embedded into $\Cal X$ by the map $t\mapsto 0\in
\Cal X(t).$ (See e.g. \cite{LT1, 1.6} for more about this deformation.)

\prop9  Assume that the tangent cone
$C_0X$ at 0 is reduced. Then, we have for $k=1,\dots,n$
$$m_k(f)=m_k(f_0),\quad\text{and }\lambda_{k+1}(f)=\lambda_{k+1}(f),$$ if
and only if for $k=2,\dots,n+1$
$$m_1(f)^k-m_k(f)=\sum_{i=2}^k\lambda_i(f)m_1(f)^{k-i}.$$ In particular, if
either condition obtains, then the smooth part of $\Cal X$ satisfies the
Whitney conditions along the parameter axis.

\pf By a result of Massey \cite{Ma, (4.7)}, for a homogeneous
polynomial of degree $d$ the following relations among the L\^e numbers
obtain:
$$\sum_{i=2}^n(d-1)^{n-i}\lambda_i(f_0)=(d-1)^{n+1}.$$ If we apply this to
$f_0$ and its restrictions to generic linear subspaces, we obtain, using
(\form.4.2), that the second relations obtain for $f=f_0.$ Therefore, the
first set of equalities implies the set of relations.

Assume now that the relations hold. We will show by induction on $k$ that
$m_k(f)=m_k(f_0)$ and $\lambda_k(f)=\lambda_k(f_0).$  First,
$m_1(f)=d-1=m_1(f_0),$ and $\lambda_1(f)=0=\lambda_1(f_0),$ as both, $f$ and
$f_0,$ are reduced by assumption. Now, assume that we have shown the equalities
for
$k-1.$ Then, the relations imply
$$\eqalign{\lambda_{k}(f)+m_{k}(f)=&m_1(f)^k-\sum_{i=2}^{k-1}\lambda_i(f)m_1
(f)^{k-i}=
m_1(f_0)^k-\sum_{i=2}^{k-1}\lambda_i(f_0)m_1(f_0)^{k-i}\cr
=&\lambda_{k}(f_0)+m_{k}(f_0).}$$
By the induction assumption and by the upper semi--continuity of Segre numbers
(\sid.5), we have $\lambda_{k}(f_0)\geq\lambda_{k}(f).$ So, it is enough
to show that
$m_{k}(f_0)\geq m_{k}(f).$  Now, the polar multiplicities $m_{k}(g_t)$
fail to be upper semi--continuous if and only if the global polar variety
$P_{k}(J_z(g),M)$ doesn't specialize, i.e
$P_{k}(J_z(g),M)(0)\neq P_{k}(J_z(g),M(0)).$ This
happens only if a component of the exceptional divisor of the blowup of
$M$ along $J_z(g)$ maps to subset of codimension $k$ in
$M(0).$

Assume that such a component $C$ exists, and let $L^k$ be a
general
$k$--plane in $(\C^{n+1},0).$ Consider the exceptional divisor $D$ of the
blowup of $L^k\times\C$ along $J_z(g).$ Then, by (\form.3), the component $C$
gives rise to a component of $D$ that maps to $0$ in
$L^k\times 0.$ But, by induction hypothesis, the L\^e numbers of
$f|_{L^k}$ equal the L\^e numbers of $f_0|_{L^k}.$ Therefore, by
(\sid.6)(i) no such component can exist.
\enddemo

\cor10 Let $X=V(f)\sub(\C^3,0)$ be a reduced hypersurface singularity. Assume
that its tangent cone $C_0X$ is reduced. Then its the deformation to the
tangent cone $\Cal X\to \C$ is Whitney--equisingular along
$0\times\C$ if and only if the following relations obtain:
$$\align
m_1(f)^2=&\lambda_2(f)+m_2(f),\\
m_1(f)^3=&\lambda_3(f)+m_1(f)\lambda_2(f).
\endalign$$

\pf This follows from the above Proposition (\wf.9) and Corollary
(\wf.6).\enddemo
\Refs
\atref{B} \by E. B\"oger
\paper Einige Bemerkungen zur Theorie der ganzalgebraischen Abh\"angikeit
von Idealen\
\jour Math. Ann. \yr 1970 \vol 185 \pages 303--308
\endref

\atref{BMM} \by Jo\"el Brian\c con, Phillipe Maisonobe, Michel Merle
\paper Localisation de syst\`ems diff\'erentiels,
stratifications de Whitney et condition de Thom\
\jour Invent. Math. \yr 1995 \vol 117 \pages 197--224
\endref

\atref{BS} \by Jo\"el Brian\c con and Jean--Paul Speder\paper Les conditions de
Whitney impliquent $\mu^*$ constant\jour Ann. Inst. Fourier \vol 26 \pages
153--163\yr 1976
\endref

\atref{Ca}
\by Roberto Callejas--Bedregal
\book Algebraic Treatment of Whitney Conditions
\publ MIT, Cambridge
\yr 1992
\bookinfo thesis
\endref

\atref{Fu}  \by William Fulton \book Intersection theory \publ
Springer-Verlag \yr 1984 \bookinfo Ergeb. Math.,3. Folge, 2. Band \endref

\atref{G1}  \by Terence Gaffney  \paper Integral Closure of Modules and Whitney
Equisingularity \jour Invent. Math. \vol 107 \yr 1992 \pages 301--322 \endref

\atref{G2}  \by Terence Gaffney  \paper Aureoles and Integral Closure of
Modules \jour preprint \yr 1992 \endref

\atref{G3}
\by Terence Gaffney
\paper Equisingularity of Plane Sections, $t_1$ Condition, and the Integral
Closure of Modules \inbook Real and Complex Singularities\ed
W.L.Marar\bookinfo Putnam Research Series in Math. 333\publ Longman\yr 1995
\endref

\atref{G4}  \by Terence Gaffney  \paper Multiplicities and equisingularity of
ICIS germs
\jour Invent. Math.
\vol 123
\yr 1996 \pages 209--220 \endref

\atref{Gas}\by L.Van Gastel \paper Excess intersections and
a correspondence principle \jour Invent.Math. \vol
103(1) \pages 197--222\yr 1991\endref

\atref{GK}
\by Terence Gaffney and Steven L. Kleiman
\paper The Principle of Specialization of Modules
\jour preprint\yr 1995
\endref

\atref{HL}
\by Jean--Pierre Georges Henry and L\^e D\~ung Tr\'ang
\paper Limites d'espaces tangents
\inbook Functions de plusieurs variable complexes II
\publ Springer Lecture Notes in Mathematics 482
\yr 1975 \bookinfo Seminaire Norguet
\pages 251--265
\endref

\atref{HMS} \by Jean--Pierre G. Henry, Michel Merle, Claude Sabbah \paper Sur
la condition de Thom stricte pour un morphisme analytique complexe \jour
Ann.Scient.\'Ec.Norm.Sup., 4. serie \vol 131 \yr 1984 \pages 227--268 \endref

\atref{KT1}
\by Steven L. Kleiman and Anders Thorup
\paper A geometric theory of the Buchsbaum--Rim multiplicity
\jour J. Algebra \vol167(1)\yr1994\pages168--231
\endref

\atref{KT2} \by Steven L. Kleiman and Anders Thorup \paper Mixed
Buchsbaum--Rim multiplicities \jour Copenhagen Univ. Preprint \vol131\yr1994
\endref

\atref{L}\by L\^e D\~ung Tr\'ang\paper Calcul du nombre de cycles
\'evanouissants d'une hypersurface complexe \jour Ann. Inst. Fourier\vol
23,4\yr 1973\pages 261--270\endref

\atref{Li}
\by Joseph Lipman
\paper Equimultiplicity, Reduction, and Blowing Up
\inbook Commutative Algebra\ed R. N. Draper  \yr 1982
\publ Marcel Dekker, New York
\bookinfo Lect. Notes Pure Appl. Math. 68
\endref

\atref{LeT}  \by Monique Lejeune-Jalabert and Bernard Teissier  \paper
Cl\^oture integrale des ideaux et equisingularit\'e, chapitre 1 \jour Publ.
Inst. Fourier  \yr 1974  
\endref

\atref{LT1} \by L\^e D\~ung Tr\'ang and Bernard Teissier  \paper Limites
d'espaces tangents en g\'eom\'etrie analytique \jour Comment. Math. Helvetici
\vol 63 \yr 1988  \pages 540--578  \endref

\atref{LT2}
\by L\^e D\~ung Tr\'ang and Bernard Teissier
\paper Vari\'et\'es polaires locales et classes de Chern des vari\'et\'es
singuli\`eres
\jour Ann. Math. \vol 114 \yr 1981  \pages 457--491
\endref

\atref{Ma} \by David Massey \book L\^e Cycles and Hypersurface Singularities
\yr 1995  \bookinfo Springer Lecture Notes in Mathematics 1615 \endref

\atref{MO} \by John Milnor and Peter Orlik\paper Isolated singularities defined
by weighted homogeneous polynomials \jour Topology \vol 9 \yr 1970 \pages
385--393
\endref

\atref{Pa}
\by Adam Parusi\'nksi\paper Limits of tangent spaces to
fibres and the $W_f$ condition\jour Duke Math. J. \vol
72,1 \yr1993\pages 99--108
\endref

\atref{R} \by David Rees \paper $A$--transforms of ideals, and a theorem on
multiplicities of ideals\jour Proc. Cambridge Phil.Soc. \vol 57\yr 1961
\pages 8--17
\endref

\atref{S} \by Pierre Samuel\paper La notion de multiplicit\'e en alg\`ebre et
en g\'eom\'etrie alg\'ebrique\jour J. Math Pures Appl. \vol 30\yr 1951
\pages 159--274\endref

\atref{T1}
\by Bernard Teissier
\paper Cycles \'evanescents, sections planes et conditions de Whitney
\jour Ast\'er\-isque \vol 7-8 \yr 1973 \pages 285--362
\endref

\atref{T2}
\by Bernard Teissier
\paper The hunting of invariants in the geometry of discriminants
\inbook Real and Complex Singularities, Oslo 1976 \publ Noordhoff \& Sijthoff
\yr 1976 \ed P. Holm \bookinfo Proceedings of the Nordic Summer School/NAVF
\pages 565--677
\endref

\atref{T3}
\by Bernard Teissier
\paper Cycles \'evanescents et resolution simultane\'e, I et II
\inbook Seminaire sur les singularit\'es des surfaces 1976--77
\publ Springer Lecture Notes in Mathematics 777
\yr 1980
\endref

\atref{T4}  \by Bernard Teissier  \paper Vari\'etes Polaires 2: Multiplicit\'es
polaires, sections planes, et conditions de Whitney \inbook Algebraic
Geometry, La R\'abida 1981 \ed J.M. Aroca  \publ Springer Lecture Notes in
Mathematics 961   \yr 1982 \bookinfo Proceedings of the International
Conference \pages 314--491 \endref

\atref{V} \by Jean--Louis Verdier \paper Stratifications de Whitney et
th\'eor\`eme de Bertini--Sard \jour Invent. Math. \vol 36 \yr 1976 \pages
295--312 \endref

\endgroup

\bigskip

 Dept.of Mathematics, Northeastern University, Boston, MA 02115, USA

 {\it E-mail address\/}: GAFF\@neu.edu, Gassler\@neu.edu
\end